\documentstyle[12pt]{article}

\font\twelveof=msym10 at 12pt
\font\sy=msxm10 at 12 pt

\def\R{\mbox{\twelveof R}}
\def\C{\mbox{\twelveof C}}
\def\N{\mbox{\twelveof N}}
\def\Z{\mbox{\twelveof Z}}
\def\sdprod{\mathrel{\mbox{\sy{\char'163}}}}
\def\id{\mbox{\rm id}}
\def\cqfd{\qquad\qquad\vrule height 4pt depth 2pt width 5pt}
\def\case#1#2{{\textstyle{#1\over #2}}}
\def\ss{\scriptstyle}
\def\arcsinh{\mathop{\rm arcsinh}\nolimits}

\newtheorem{definition}{Definition}[section]
\newtheorem{proposition}[definition]{Proposition}
\newtheorem{corollary}[definition]{Corollary}

\title{\hspace{9cm}{\normalsize ULB/229/CQ/96/5}\\
\vspace{1.5cm}
Coloured quantum universal enveloping algebras}
\author{C. Quesne\thanks{Directeur de recherches FNRS; E-mail:
cquesne@ulb.ac.be} \\
{\small \sl Physique Nucl\'eaire Th\'eorique et Physique Math\'ematique, 
Universit\'e Libre de Bruxelles,} \\
{\small \sl Campus de la Plaine CP229, Boulevard~du Triomphe, B-1050 Brussels,
Belgium}}
\date{ }
\begin{document}
\baselineskip=22pt plus 1pt minus 1pt
\maketitle

\begin{abstract}
We define some new algebraic structures, termed coloured Hopf algebras, by
combining the coalgebra structures and antipodes of a standard Hopf algebra
set~$\cal H$, corresponding to some parameter set~$\cal Q$, with the
transformations of an algebra isomorphism group~$\cal G$, herein called colour
group. Such transformations are labelled by some colour parameters, taking values
in a colour set~$\cal C$. We show that various classes of Hopf algebras, such as
almost cocommutative, coboundary, quasitriangular, and triangular ones, can be
extended into corresponding coloured algebraic structures, and that coloured
quasitriangular Hopf algebras, in particular, are characterized by the existence of a
coloured universal $\cal R$-matrix, satisfying the coloured Yang-Baxter
equation. The present definitions extend those previously introduced by Ohtsuki,
which correspond to some substructures in those cases where the colour group is
abelian. We apply the new concepts to construct coloured quantum universal
enveloping algebras of both semisimple and nonsemisimple Lie algebras, 
considering several examples with fixed or varying parameters. As a
by-product, some of the matrix representations of coloured universal
$\cal R$-matrices, derived in the present paper, provide new solutions of the
coloured Yang-Baxter equation, which might be of interest in the context of
integrable models.
\end{abstract}

\vspace{0.5cm}

\hspace*{0.3cm}

\hspace*{0.3cm}
Running title: Coloured quantum algebras

\hspace*{0.3cm}
To appear in J. Math. Phys.
\newpage
%
%
\section{INTRODUCTION} \label{sec:intro}
Since its introduction, the parametrized (quantum) Yang-Baxter equation
(YBE)~\cite{yang} plays a crucial role in nonlinear integrable systems in physics,
such as exactly solvable statistical mechanics models and low-dimensional
integrable field theories~\cite{kulish}. Its constant form is also important in knot
theory, where it is connected with braid groups~\cite{jones}.\par
%
%
In addition, the~YBE has inspired the development of quantum groups and quantum
algebras~\cite{sklyanin}, essentially appearing in the literature in two different
forms.\par
%
%
In the Faddeev-Reshetikhin-Takhtajan (FRT) formulation~\cite{faddeev}, to any
invertible matrix solution~$R$ of the constant YBE, one associates two
bialgebras~$A(R)$ and $U(R)$ that under certain conditions can lead to two dual
Hopf algebras. For both of them, the constant YBE is a sufficient condition for
associativity.\par 
%
%
In the Drinfeld and Jimbo (DJ) approach~\cite{drinfeld}, one considers
one-parameter deformations~$U_q(g)$ of the universal enveloping algebras (or
quantum universal enveloping algebras (QUEA's)) of simple Lie algebras~$g$. Such
quantizations are quasitriangular Hopf algebras, that is there exists a universal
${\cal R}$-matrix, which is an invertible element of $U_q(g) \otimes U_q(g)$, and
which among other properties, satisfies the YBE.\par
%
%
The DJ approach has been completed for nonsemisimple Lie algebras by applying
various procedures, such as contractions of QUEA's of simple Lie
algebras~\cite{celeghini90, celeghini91}, and by introducing multiparametric
deformations~\cite{schirrmacher, dobrev}.\par
%
%
A link has been established between the FRT and DJ formulations by considering
for~$R$ the matrix representing the operator~${\cal R}$ in the fundamental
representation of~$g$ (see e.g.\ Ref.~\cite{burroughs}). In such a case, the
generators of~$U(R)$ can indeed be expressed in terms of those of~$U_q(g)$. The
algebra~$A(R)$, associated with~$U(R)$, is then related to the quantized functions
on the corresponding Lie group~$G$, or quantum group, $Fun_q(G) = G_q$.\par
%
%
In recent years, some integrable models with nonadditive-type solutions
$R^{\lambda,\mu} \ne R(\lambda - \mu)$ of the~YBE have been
discovered~\cite{bazhanov, hlavaty}. The corresponding~YBE
\begin{equation}
  R^{\lambda,\mu}_{12} R^{\lambda,\nu}_{13} R^{\mu,\nu}_{23} = R^{\mu,\nu}_{23}
  R^{\lambda,\nu}_{13} R^{\lambda,\mu}_{12}   \label{eq:colYBE}
\end{equation}
is referred to in the literature as the `coloured'~YBE, the nonadditive (in general
multicomponent) spectral parameters $\lambda$, $\mu$,~$\nu$ being considered
as `colour' indices~\cite{footnote1}.\par
%
%
Constructing solutions of Eq.~(\ref{eq:colYBE}) by starting from some quantum
algebra has then become a topic of active research. Various approaches have been
used for such a purpose~\cite{akutsu}--\cite{bonatsos}. Among them, one should
mention a recent work of Bonatsos {\it et al\/}~\cite{bonatsos} on a nonlinear
deformation
${\cal A}^+_q(1)$ of $su(2)$, distinct from the DJ one, wherein the colour parameter
is related with an involutive automorphism of the algebra, and serves to distinguish
between the irreducible representations with the same dimension.\par
%
%
Extending the definitions of quantum groups and quantum algebras by connecting
them to coloured $R$-matrices, instead of ordinary ones, is an interesting
problem, which so far has not received much attention in the literature.\par
%
%
The generalization of the FRT approach has been discussed by Kundu and
Basu-Mallick~\cite{basu94, kundu94, basu95} for some quantizations of $U(gl(2))$
and $Gl(2)$. Such coloured extensions are characterized by generalized algebraic
structures, but coalgebra structures identical with those of standard $A(R)$ and
$U(R)$ algebras.\par
%
%
In the context of knot theory, Ohtsuki~\cite{ohtsuki} has introduced coloured
quasitriangular Hopf algebras, which are characterized by the existence of a
coloured universal $\cal R$-matrix, and he has applied his formalism to coloured
representations of $U_q(sl(2))$ for $q$ a root of unity. A rather similar, but
nevertheless distinct generalization has been independently considered by Bonatsos
{\it et al\/}~\cite{bonatsos} for the above-mentioned ${\cal A}^+_q(1)$ algebra,
which has been endowed with a two-colour quasitriangular Hopf structure.\par
%
%
In the present paper, we extend the DJ formulation of QUEA's to coloured ones by
elaborating on the results of Bonatsos {\it et~al}~\cite{bonatsos}. For
such a purpose, in the next section we define coloured Hopf algebras in a way that
generalizes Ohtsuki's first attempt. In Secs.~\ref{sec:fixed} and~\ref{sec:varying},
we demonstrate the present definition generality and usefulness by reviewing
various examples of coloured QUEA's, then summarize and comment on prospects in
Sec.~\ref{sec:conclusion}.\par
%
%
\section{COLOURED HOPF ALGEBRAS}    \label{sec:colHopf}
\setcounter{equation}{0}
Let $\left({\cal H}_q, +, m_q, \iota_q, \Delta_q, \epsilon_q, S_q;k\right)$ (or in
short ${\cal H}_q$) be a Hopf algebra over some field~$k$ ($= \C$ or \R), depending
upon some parameters~$q$. Here $m_q: {\cal H}_q \otimes {\cal H}_q \to {\cal
H}_q$, $\iota_q: k \to {\cal H}_q$, $\Delta_q: {\cal H}_q \to {\cal H}_q \otimes
{\cal H}_q$, $\epsilon_q: {\cal H}_q \to k$, and $S_q: {\cal H}_q \to {\cal H}_q$
denote the multiplication, unit, comultiplication, counit, and antipode maps
respectively~\cite{majid}. Whenever $q$ runs over some set~$\cal Q$, called
{\it parameter set\/}, we obtain a set of Hopf algebras ${\cal H} = \{\,{\cal H}_q
\mid q
\in {\cal Q}\,\}$. In the examples given in Secs.~\ref{sec:fixed}
and~\ref{sec:varying}, we shall distinguish between two cases, according to
whether $\cal Q$ contains a single element (fixed-parameter case) or more than one
element (varying-parameter case).\par
%
%
Let us assume that there exists a set of one-to-one linear maps ${\cal G} = \{\,
\sigma^{\nu}: {\cal H}_q \to {\cal H}_{q^{\nu}} \mid q, q^{\nu} \in {\cal Q}, \nu \in
{\cal C}\,\}$, defined for any ${\cal H}_q \in {\cal H}$. They are labelled by some
parameters~$\nu$, called {\it colour parameters\/}, taking values in some
set~$\cal C$, called {\it colour set\/}. The latter may be finite, countably infinite,
or uncountably infinite. Two conditions are imposed on the $\sigma^{\nu}$'s:

\begin{itemize}
\item[(i)] Every $\sigma^{\nu}$ is an algebra isomorphism, i.e.,
\begin{equation}
  \sigma^{\nu} \circ m_q = m_{q^{\nu}} \circ \left(\sigma^{\nu} \otimes
  \sigma^{\nu}\right), \qquad \sigma^{\nu} \circ \iota_q = \iota_{q^{\nu}};
  \label{eq:iso} 
\end{equation}
\item[(ii)] $\cal G$ is a group (called {\it colour group\/}) with respect to the
composition of maps, i.e.,
\begin{eqnarray}
  \forall \nu, \nu' \in {\cal C}, \exists\, \nu'' \in {\cal C}: \sigma^{\nu''} & =
           &\sigma^{\nu'} \circ \sigma^{\nu}: {\cal H}_q \to {\cal H}_{q^{\nu''}} = 
          {\cal H}_{q^{\nu,\nu'}},  \label{eq:compo} \\
  \exists\, \nu^0 \in {\cal C}: \sigma^{\nu^0} & = & \id: {\cal H}_q \to 
          {\cal H}_{q^{\nu^0}} = {\cal H}_q, \\
  \forall \nu \in {\cal C}, \exists\, \nu' \in {\cal C}: \sigma^{\nu'} & =
           &\sigma_{\nu} \equiv \left(\sigma^{\nu}\right)^{-1}: {\cal H}_{q^{\nu}} \to
          {\cal H}_{q}   \label{eq:inverse}.
\end{eqnarray}
In Eqs.~(\ref{eq:compo}) and (\ref{eq:inverse}), $\nu''$ and $\nu'$ will be denoted by
$\nu' \circ \nu$ and $\nu^i$, respectively.
\end{itemize}
\par
%
%
$\cal H$, $\cal C$, and $\cal G$ can be combined into
%
\begin{definition}   \label{def-colmaps}
The maps $\Delta^{\lambda,\mu}_{q,\nu}: {\cal H}_{q^{\nu}} \to {\cal
H}_{q^{\lambda}} \otimes {\cal H}_{q^{\mu}}$, $\epsilon_{q,\nu}: {\cal H}_{q^{\nu}}
\to k$, and $S^{\mu}_{q,\nu}: {\cal H}_{q^{\nu}} \to {\cal H}_{q^{\mu}}$, defined by
\begin{equation}
  \Delta^{\lambda,\mu}_{q,\nu} \equiv \left(\sigma^{\lambda} \otimes
  \sigma^{\mu}\right) \circ \Delta_q \circ \sigma_{\nu}, \qquad
  \epsilon_{q,\nu} \equiv \epsilon_q \circ \sigma_{\nu}, \qquad
  S^{\mu}_{q,\nu} \equiv \sigma^{\mu} \circ S_q \circ \sigma_{\nu},
  \label{eq:colmaps} 
\end{equation}
for any $q \in {\cal Q}$, and any $\lambda$, $\mu$, $\nu \in {\cal C}$, are called
coloured comultiplication, counit, and antipode respectively. 
\end{definition}
\par
%
%
It is easy to prove the following proposition:
%
\begin{proposition}   \label{prop-GenHopf}
The coloured comultiplication, counit, and antipode maps, defined in
Eq.~(\ref{eq:colmaps}), transform under the colour group~$\cal G$ as
\begin{eqnarray}
  \left(\sigma^{\lambda}_{\alpha} \otimes \sigma^{\mu}_{\beta}\right) \circ 
           \Delta^{\alpha,\beta}_{q,\nu} & = & \Delta^{\lambda,\mu}_{q,\nu} = 
           \Delta^{\lambda,\mu}_{q,\gamma} \circ \sigma^{\gamma}_{\nu}, 
           \nonumber \\
  \epsilon_{q,\alpha} \circ \sigma^{\alpha}_{\nu} & = & \epsilon_{q,\nu},
           \nonumber \\
  \sigma^{\mu}_{\alpha} \circ S^{\alpha}_{q,\nu} & = & S^{\mu}_{q,\nu} = 
           S^{\mu}_{q,\beta} \circ \sigma^{\beta}_{\nu},
\end{eqnarray}
and satisfy generalized coassociativity, counit, and antipode axioms
\begin{eqnarray}
  \left(\Delta^{\alpha,\beta}_{q,\lambda} \otimes \sigma^{\gamma}_{\mu}\right)
           \circ \Delta^{\lambda,\mu}_{q,\nu} & = & \left(\sigma^{\alpha}_{\lambda'}
           \otimes \Delta^{\beta,\gamma}_{q,\mu'}\right) \circ 
           \Delta^{\lambda',\mu'}_{q,\nu}, \nonumber \\
  \left(\epsilon_{q,\lambda} \otimes \sigma^{\alpha}_{\mu}\right) \circ
           \Delta^{\lambda,\mu}_{q,\nu} & = & \left(\sigma^{\alpha}_{\lambda'}
           \otimes \epsilon_{q,\mu'}\right) \circ \Delta^{\lambda',\mu'}_{q,\nu} =
           \sigma^{\alpha}_{\nu}, \nonumber \\
  m_{q^{\alpha}} \circ \left(S^{\alpha}_{q,\lambda} \otimes \sigma^{\alpha}_{\mu}
           \right) \circ \Delta^{\lambda,\mu}_{q,\nu} & = & m_{q^{\alpha}} \circ \left(
           \sigma^{\alpha}_{\lambda'} \otimes S^{\alpha}_{q,\mu'} \right) \circ
           \Delta^{\lambda',\mu'}_{q,\nu} = \iota_{q^{\alpha}} \circ \epsilon_{q,\nu},
           \label{eq:colantipode}
\end{eqnarray}
as well as generalized bialgebra axioms
\begin{eqnarray}
  \Delta^{\lambda,\mu}_{q,\nu} \circ m_{q^{\nu}} & = & \left(m_{q^{\lambda}}
           \otimes m_{q^{\mu}}\right) \circ (\id \otimes \tau \otimes \id) \circ
           \left(\Delta^{\lambda,\mu}_{q,\nu} \otimes \Delta^{\lambda,\mu}_{q,\nu}
           \right), \nonumber \\
  \Delta^{\lambda,\mu}_{q,\nu} \circ \iota_{q^{\nu}} & = & \iota_{q^{\lambda}}
           \otimes \iota_{q^{\mu}}, \nonumber \\
  \epsilon_{q,\nu} \circ m_{q^{\nu}} & = & \epsilon_{q,\nu} \otimes \epsilon_{q,\nu},
           \nonumber \\
  \epsilon_{q,\nu} \circ \iota_{q^{\nu}} & = & 1_k.
\end{eqnarray}
Here $\sigma^{\lambda}_{\mu}$ is the element of~$\cal G$ defined by
\begin{equation}
  \sigma^{\lambda}_{\mu} \equiv \sigma^{\lambda} \circ \sigma_{\mu},
\end{equation}
$\tau$ is the twist map, i.e., $\tau(a \otimes b) = b \otimes a$, $1_k$ denotes the
unit of~$k$, and no summation is implied over repeated indices.
\end{proposition}
%
{\it Proof.\/} The various results are obtained by combining standard Hopf algebra
axioms~\cite{majid} with Definition~\ref{def-colmaps}, and
Eqs.~(\ref{eq:iso})--(\ref{eq:inverse}). Consider for instance the first equation
in~(\ref{eq:colantipode}). The map on the left-hand side, $\left(\Delta^{\alpha,\beta}
_{q,\lambda} \otimes \sigma^{\gamma}_{\mu}\right) \circ
\Delta^{\lambda,\mu}_{q,\nu}: {\cal H}_{q^{\nu}} \to {\cal H}_{q^{\lambda}} \otimes
{\cal H}_{q^{\mu}} \to {\cal H}_{q^{\alpha}} \otimes {\cal H}_{q^{\beta}} \otimes
{\cal H}_{q^{\gamma}}$, can be proved to be identical with that on the right-hand
one, $\left(\sigma^{\alpha}_{\lambda'} \otimes
\Delta^{\beta,\gamma}_{q,\mu'}\right) \circ \Delta^{\lambda',\mu'}_{q,\nu}: 
{\cal H}_{q^{\nu}} \to {\cal H}_{q^{\lambda'}} \otimes {\cal H}_{q^{\mu'}} \to {\cal
H}_{q^{\alpha}} \otimes {\cal H}_{q^{\beta}} \otimes {\cal H}_{q^{\gamma}}$, as
follows:
\begin{eqnarray*}
  \lefteqn{\left(\Delta^{\alpha,\beta}_{q,\lambda} \otimes
          \sigma^{\gamma}_{\mu}\right) \circ \Delta^{\lambda,\mu}_{q,\nu} =
           \left(\left(\left(\sigma^{\alpha} \otimes \sigma^{\beta}\right) \circ
          \Delta_q \circ \sigma_{\lambda}\right) \otimes \left(\sigma^{\gamma}
          \circ \sigma_{\mu}\right)\right) \circ \left(\left(\sigma^{\lambda} \otimes
          \sigma^{\mu}\right) \circ \Delta_q \circ \sigma_{\nu}\right)} \\
  & = & \left(\left(\left(\sigma^{\alpha} \otimes \sigma^{\beta}\right) \circ
          \Delta_q\right) \otimes \sigma^{\gamma}\right) \circ \Delta_q \circ
          \sigma_{\nu} = \left(\sigma^{\alpha} \otimes \sigma^{\beta} \otimes
          \sigma^{\gamma}\right) \circ \left(\Delta_q \otimes \id\right)  \circ
          \Delta_q \circ \sigma_{\nu} \\
  & = & \left(\sigma^{\alpha} \otimes \sigma^{\beta} \otimes
          \sigma^{\gamma}\right) \circ \left(\id \otimes \Delta_q\right) \circ
          \Delta_q \circ \sigma_{\nu} = \left(\sigma^{\alpha} \otimes \left(\left(
          \sigma^{\beta} \otimes \sigma^{\gamma}\right) \circ \Delta_q\right)\right)
          \circ \Delta_q \circ \sigma_{\nu} \\
  & = & \left(\left(\sigma^{\alpha} \circ \sigma_{\lambda'}\right) \otimes
          \left(\left(\sigma^{\beta} \otimes \sigma^{\gamma}\right) \circ \Delta_q
          \circ \sigma_{\mu'}\right)\right) \circ \left(\sigma^{\lambda'} \otimes
          \sigma^{\mu'}\right) \circ \Delta_q \circ \sigma_{\nu} \\
  & = & \left(\sigma^{\alpha}_{\lambda'} \otimes
          \Delta^{\beta,\gamma}_{q,\mu'}\right) \circ
          \Delta^{\lambda',\mu'}_{q,\nu}. \cqfd    
\end{eqnarray*}
\par
%
%
{}From Proposition~\ref{prop-GenHopf}, it is straightforward to obtain
%
\begin{corollary}   \label{corol-newHopf}
If Eqs.~(\ref{eq:iso})--(\ref{eq:inverse}) are satisfied, then for any $q \in \cal Q$,
any $\nu \in \cal C$, and $q_{\nu} \equiv q^{\nu^i}$, $\left({\cal H}_q, +, m_q,
\iota_q, \Delta^{\nu,\nu}_{q_{\nu},\nu}, \epsilon_{q_{\nu},\nu},
S^{\nu}_{q_{\nu},\nu}; k\right)$ is a Hopf algebra over $k$ with comultiplication
$\Delta^{\nu,\nu}_{q_{\nu},\nu}$, counit $\epsilon_{q_{\nu},\nu}$, and antipode
$S^{\nu}_{q_{\nu},\nu}$, defined by particularizing Eq.~(\ref{eq:colmaps}).
\end{corollary}
%
{\it Remark.\/} In particular, for $\nu = \nu^0$, we get back the original Hopf
structure of~${\cal H}_q$.\par
%
%
Generalizing the result contained in Corollary~\ref{corol-newHopf}, we are led to
introduce
%
\begin{definition}
A set of Hopf algebras~$\cal H$, endowed with coloured comultiplication, counit,
and antipode maps $\Delta^{\lambda,\mu}_{q,\nu}$, $\epsilon_{q,\nu}$,
$S^{\mu}_{q,\nu}$, as defined in (\ref{eq:colmaps}), is called coloured Hopf algebra,
and denoted by any one of the symbols $\left({\cal H}_q, +, m_q, \iota_q,
\Delta^{\lambda,\mu}_{q,\nu}, \epsilon_{q,\nu}, S^{\mu}_{q,\nu}; k, {\cal Q}, {\cal
C}, {\cal G}\right)$, $\left({\cal H}, {\cal C}, {\cal G}\right)$, or ${\cal H}^c$.
\end{definition}
\par
%
%
As in standard Hopf algebras, the coloured antipode~$S^{\mu}_{q,\nu}$ satisfies
some additional properties.
%
\begin{proposition}
The coloured antipode~$S^{\mu}_{q,\nu}$ of a coloured Hopf algebra~${\cal H}^c$
fulfils the relations
\begin{eqnarray}
  S^{\mu}_{q,\nu} \circ m_{q^{\nu}} & = & m_{q^{\mu}} \circ \tau \circ \left(
           S^{\mu}_{q,\nu} \otimes S^{\mu}_{q,\nu}\right), \qquad S^{\mu}_{q,\nu}
           \circ \iota_{q^{\nu}} = \iota_{q^{\mu}},  \label{eq:S1} \\
  \left(S^{\alpha}_{q,\lambda} \otimes S^{\beta}_{q,\mu}\right) \circ
           \Delta^{\lambda,\mu}_{q,\nu} & = & \tau \circ
           \Delta^{\beta,\alpha}_{q,\gamma} \circ S^{\gamma}_{q,\nu}, \qquad
           \epsilon_{q,\mu} \circ S^{\mu}_{q,\nu} = \epsilon_{q,\nu}.  \label{eq:S2}
\end{eqnarray}
\end{proposition}
%
{\it Proof.\/} Eq.~(\ref{eq:S1}) (resp.~(\ref{eq:S2})) is obtained by combining 
Eqs.~(\ref{eq:iso})--(\ref{eq:colmaps}) with the first (resp.~second) line in the
following equation
\begin{eqnarray*}
  S_q \circ m_q & = & m_q \circ \tau \circ \left(S_q \otimes S_q\right), \qquad
           S_q \circ \iota_q = \iota_q, \\
\left(S_q \otimes S_q\right) \circ \Delta_q & = & \tau \circ \Delta_q \circ S_q,
           \qquad \epsilon_q \circ S_q = \epsilon_q,
\end{eqnarray*}
expressing the fact that $S_q$ is an algebra (resp.~coalgebra)
antiautomorphism. \cqfd\par
%
%
Let us now assume that the members of the Hopf algebra set~$\cal H$ are almost
cocommutative Hopf algebras~\cite{majid}, i.e., for any $q \in \cal Q$ there exists
an invertible element ${\cal R}_q \in {\cal H}_q \otimes {\cal H}_q$ (completed
tensor product), such that
\begin{equation}
  \tau \circ \Delta_q(a) = {\cal R}_q \Delta_q(a) {\cal R}_q^{-1}
\end{equation}
for any $a \in {\cal H}_q$.\par
%
%
We may then introduce
%
\begin{definition}    \label{def-colR}
Let ${\cal R}^c$ denote the set of elements ${\cal R}^{\lambda,\mu}_q \in {\cal
H}_{q^{\lambda}} \otimes {\cal H}_{q^{\mu}}$, defined by
\begin{equation}
  {\cal R}^{\lambda,\mu}_q \equiv \left(\sigma^{\lambda} \otimes
  \sigma^{\mu}\right) \left({\cal R}_q\right),    \label{eq:colR}  
\end{equation}
where $q$ runs over~$\cal Q$, and $\lambda$, $\mu$ over~$\cal C$.
\end{definition}
%
The following result can be easily obtained:
%
\begin{proposition}
If the Hopf algebras~${\cal H}_q$ of~$\cal H$ are almost cocommutative, then
${\cal R}^{\lambda,\mu}_q$, as defined in (\ref{eq:colR}), is invertible with
$\left({\cal R}^{\lambda,\mu}_q\right)^{-1}$ given by
\begin{equation}
  \left({\cal R}^{\lambda,\mu}_q\right)^{-1} = \left(\sigma^{\lambda} \otimes
  \sigma^{\mu}\right) \left({\cal R}_q^{-1}\right),   \label{eq:alcocom1}
\end{equation}
and
\begin{equation}
  \tau \circ \Delta^{\mu,\lambda}_{q,\nu}(a) = {\cal R}^{\lambda,\mu}_q
  \Delta^{\lambda,\mu}_{q,\nu}(a) \left({\cal R}^{\lambda,\mu}_q\right)^{-1}
  \label{eq:alcocom2}
\end{equation}
for any $a \in {\cal H}_{q^{\nu}}$. If in addition, the almost cocommutative Hopf
algebras $\left({\cal H}_q, {\cal R}_q\right)$ are (i) coboundary, (ii)
quasitriangular, or (iii) triangular, then ${\cal R}^{\lambda,\mu}_q$ also satisfies
the relations
\begin{itemize}
\item[(i)]
\begin{eqnarray}
  {\cal R}^{\alpha,\beta}_{q,12} \left(\Delta^{\alpha,\beta}_{q,\lambda} \otimes
             \sigma^{\gamma}_{\mu}\right) \left({\cal R}^{\lambda,\mu}_q\right) & = &
             {\cal R}^{\beta,\gamma}_{q,23} \left(\sigma^{\alpha}_{\lambda'} \otimes
             \Delta^{\beta,\gamma}_{q,\mu'}\right) \left({\cal R}^{\lambda',\mu'}_q
             \right), \nonumber \\
  {\cal R}^{\lambda,\mu}_{q,21} & \equiv & \tau\left({\cal
              R}^{\mu,\lambda}_q\right) = \left({\cal R}^{\lambda,\mu}_q\right)^{-1},
              \nonumber \\
  \left(\epsilon_{q,\lambda} \otimes \epsilon_{q,\mu}\right) \left({\cal
             R}^{\lambda,\mu}_q\right) & = & 1_k,  \label{eq:coboundary}
\end{eqnarray}
\item[(ii)] 
\begin{eqnarray}
  \left(\Delta^{\alpha,\beta}_{q,\lambda} \otimes \sigma^{\gamma}_{\mu}\right)
            \left({\cal R}^{\lambda,\mu}_q\right) & = & {\cal R}^{\alpha,\gamma}_{q,13}
            {\cal R}^{\beta,\gamma}_{q,23}, \nonumber \\
  \left(\sigma^{\alpha}_{\lambda} \otimes \Delta^{\beta,\gamma}_{q,\mu}\right)
            \left({\cal R}^{\lambda,\mu}_q\right) & = & {\cal R}^{\alpha,\gamma}_{q,13}
            {\cal R}^{\alpha,\beta}_{q,12},    \label{eq:quasi}
\end{eqnarray}
\item[(iii)] 
\begin{eqnarray}
  \left(\Delta^{\alpha,\beta}_{q,\lambda} \otimes \sigma^{\gamma}_{\mu}\right)
            \left({\cal R}^{\lambda,\mu}_q\right) & = & {\cal R}^{\alpha,\gamma}_{q,13}
            {\cal R}^{\beta,\gamma}_{q,23}, \nonumber \\
  \left(\sigma^{\alpha}_{\lambda} \otimes \Delta^{\beta,\gamma}_{q,\mu}\right)
            \left({\cal R}^{\lambda,\mu}_q\right) & = & {\cal R}^{\alpha,\gamma}_{q,13}
            {\cal R}^{\alpha,\beta}_{q,12}, \nonumber \\
  {\cal R}^{\lambda,\mu}_{q,21} & = & \left({\cal R}^{\lambda,\mu}_q\right)^{-1},
            \label{eq:triangular}
\end{eqnarray}
\end{itemize}
respectively.
\end{proposition}
%
Hence we have
%
\begin{definition}   \label{def-quasi}
A coloured, almost cocommutative Hopf algebra is a pair $\left({\cal H}^c, {\cal
R}^c\right)$, where ${\cal H}^c$ is a coloured Hopf algebra, ${\cal R}^c = \{\, {\cal
R}^{\lambda,\mu}_q \mid q \in {\cal Q}, \lambda, \mu \in {\cal C}\,\}$, and ${\cal
R}^{\lambda,\mu}_q$, defined in (\ref{eq:colR}), satisfies Eqs.~(\ref{eq:alcocom1})
and~(\ref{eq:alcocom2}). A coloured, almost cocommutative Hopf algebra 
$\left({\cal H}^c, {\cal R}^c\right)$ is said to be coboundary, quasitriangular, or
triangular if ${\cal R}^{\lambda,\mu}_q$ satisfies Eq.~(\ref{eq:coboundary}),
(\ref{eq:quasi}), or~(\ref{eq:triangular}), respectively. In the case of a coloured
quasitriangular Hopf algebra, the set~${\cal R}^c$ is called the coloured universal
$\cal R$-matrix of $\left({\cal H}^c, {\cal R}^c\right)$. 
\end{definition}
\par
%
%
The terminology used for~${\cal R}^c$ in Definition~\ref{def-quasi} is justified by
the following proposition:
%
\begin{proposition}
Let $\left({\cal H}^c, {\cal R}^c\right)$ be a coloured quasitriangular Hopf algebra.
Then
\begin{eqnarray}
  {\cal R}^{\lambda,\mu}_{q,12} {\cal R}^{\lambda,\nu}_{q,13} 
            {\cal R}^{\mu,\nu}_{q,23} & = & {\cal R}^{\mu,\nu}_{q,23}
            {\cal R}^{\lambda,\nu}_{q,13} {\cal R}^{\lambda,\mu}_{q,12}, \nonumber \\
  \left(\epsilon_{q,\lambda} \otimes \sigma^{\alpha}_{\mu}\right) \left({\cal
            R}^{\lambda,\mu}_q\right) & = & \left(\sigma^{\alpha}_{\lambda'} \otimes
            \epsilon_{q,\mu'}\right) \left({\cal R}^{\lambda',\mu'}_q\right) =
           1_{q^{\alpha}}, \nonumber \\
  \left(S^{\alpha}_{q,\lambda} \otimes \sigma^{\beta}_{\mu}\right) \left({\cal
            R}^{\lambda,\mu}_q\right) & = & \left(\sigma^{\alpha}_{\lambda'} \otimes
            \left(S^{\mu'}_{q,\beta}\right)^{-1}\right) \left({\cal
            R}^{\lambda',\mu'}_q\right) = \left({\cal
            R}^{\alpha,\beta}_q\right)^{-1},   \label{eq:propquasi}  
\end{eqnarray}
where $1_{q^{\alpha}}$ denotes the unit element of~${\cal H}_{q^{\alpha}}$, and
$\left(S^{\mu}_{q,\nu}\right)^{-1}: {\cal H}_{q^{\mu}} \to {\cal H}_{q^{\nu}}$ is
given by $\left(S^{\mu}_{q,\nu}\right)^{-1} = \sigma^{\nu} \circ S_q^{-1} \circ
\sigma_{\mu}$.
\end{proposition}
%
{\it Remarks.\/} (1) The first equation in~(\ref{eq:propquasi}) shows that the
elements of the coloured universal $\cal R$-matrix satisfy the coloured YBE, as
given in~(\ref{eq:colYBE}). (2) Following common use for standard Hopf
algebras, we shall also call ${\cal R}^c$ coloured universal $\cal R$-matrix if its
elements satisfy Eq.~(\ref{eq:alcocom2}) and the coloured YBE.\par
%
%
In those cases where the colour group~$\cal G$ is abelian, one can always
transform the colour parameters so as to make them additive. Let therefore
$\nu(p)$ be such that $\nu'(p') \circ \nu(p) = (\nu' \circ \nu)(p+p')$, $\nu(0) =
\nu^0$, $\nu(-p) = \nu^i(p)$, and let denote ${\cal H}_{q^{\nu(p)}}$ by ${\cal A}_p$.
The coloured comultiplication, counit, antipode, and universal $\cal R$-matrix,
introduced in Definitions~\ref{def-colmaps} and~\ref{def-colR}, can then be
written as $\Delta^{\lambda(p_1),\mu(p_2)}_{q,\nu(p_3)}$, $\epsilon_{q,\nu(p)}$,
$S^{\mu(p_1)}_{q,\nu(p_2)}$, and ${\cal R}^{\lambda(p_1),\mu(p_2)}_q$,
respectively. By specializing the results obtained in (\ref{eq:colantipode}),
(\ref{eq:alcocom2}), and~(\ref{eq:quasi}), we obtain
%
\begin{proposition}
If $\left({\cal H}^c, {\cal R}^c\right)$ is a coloured quasitriangular Hopf algebra
with an abelian colour group~$\cal G$, then the maps $\Delta_{p_1p_2} \equiv 
\Delta^{\lambda(p_1),\mu(p_2)}_{q,(\lambda \circ \mu)(p_1+p_2)}: {\cal
A}_{p_1+p_2} \to {\cal A}_{p_1} \otimes {\cal A}_{p_2}$, $\epsilon \equiv
\epsilon_{q,\nu(0)}: {\cal A}_0 \to k$, $S_p \equiv S^{\nu(-p)}_{q,\nu(p)}: {\cal A}_p
\to {\cal A}_{-p}$, and the invertible elements ${\cal R}_{p_1p_2} \equiv {\cal
R}^{\lambda(p_1),\mu(p_2)}_q$ of ${\cal A}_{p_1} \otimes {\cal A}_{p_2}$ satisfy
the defining relations of an Ohtsuki's coloured quasitriangular Hopf
algebra~\cite{ohtsuki}, i.e.,
\begin{eqnarray}
  \left(\Delta_{p_1p_2} \otimes \id\right) \circ \Delta_{p_1+p_2,p_3} & = &
           \left(\id \otimes \Delta_{p_2p_3}\right) \circ \Delta_{p_1,p_2+p_3},
           \nonumber \\
  \left(\epsilon \otimes \id\right) \circ \Delta_{0,p} & = & \left(\id \otimes
           \epsilon\right) \circ \Delta_{p,0} = \id, \nonumber \\
  m_p \circ \left(S_{-p} \otimes \id\right) \circ \Delta_{-p,p} & = & m_p \circ
           \left(\id \otimes S_{-p}\right) \circ \Delta_{p,-p} = \iota_p \circ \epsilon,
           \nonumber \\
  \tau \circ \Delta_{p_2p_1}(a) & = & {\cal R}_{p_1p_2} \Delta_{p_1p_2}(a)
           {\cal R}_{p_1p_2}^{-1}, \nonumber \\
  \left(\Delta_{p_1p_2} \otimes \id \right) \left({\cal R}_{p_1+p_2,p_3}\right) 
           & = & {\cal R}_{p_1p_3,13} {\cal R}_{p_2p_3,23}, \nonumber \\
  \left(\id \otimes \Delta_{p_2p_3}\right) \left({\cal R}_{p_1,p_2+p_3}\right) 
           & = & {\cal R}_{p_1p_3,13} {\cal R}_{p_1p_2,12},
\end{eqnarray}
with $m_p \equiv m_{q^{\nu(p)}}$, $\iota_p \equiv \iota_{q^{\nu(p)}}$, and $a \in
{\cal A}_{p_1+p_2}$. Hence the latter is a substructure of $\left({\cal H}^c, {\cal
R}^c\right)$.
\end{proposition}
%
{\it Remark.\/} As opposed to Ohtsuki's coloured Hopf algebras, those considered in
the present paper are also valid for nonabelian colour groups. Such a generalization
is significant as it will be shown in the next two sections that coloured Hopf
algebras with such colour groups can indeed be constructed.\par
%
%
\section{EXAMPLES OF COLOURED QUEA'S WITH FIX\-ED PARAMETERS} 
\label{sec:fixed}
\setcounter{equation}{0}
In the present section, we construct various examples of coloured
quasitriangular Hopf algebras, for which the underlying Hopf algebra set~$\cal H$
reduces to a single element~${\cal H}_q$ (hence ${\cal Q} = \{q\}$ and $q^{\nu} =
q$), which is some QUEA $U_q(g)$.\par
%
%
\subsection{The standard quantum algebra $U_q(sl(2))$}  
\label{subsec:sl(2)}
We begin by considering the simplest case of QUEA, namely the standard
DJ~deformation of $U(sl(2))$~\cite{drinfeld}, i.e., $U_q(sl(2))$ where $k = \C$ and
$q = \exp(\eta) \in \C \! \setminus \! \{0\}$, whose universal $\cal R$-matrix was
obtained in Ref.~\cite{kirillov}. Although this example might look over-simple, it
nevertheless serves three important purposes: to illustrate the fact that any
QUEA can be easily transformed into a coloured one, to show that this can be
achieved in various ways, and to demonstrate that some of them may involve a
nonabelian colour group.\par
%
%
\subsubsection{The colour group ${\cal G} = S_2$}   \label{subsubsec:S_2}
The quantum algebra $U_q(sl(2))$ is generated by three operators
$J_3$,~$J_{\pm}$, satisfying the commutation relations
\begin{equation}
  \left[J_3, J_{\pm}\right] = \pm J_{\pm}, \qquad \left[J_+, J_-\right] = 
  \left[2J_3\right]_q = \frac{\sinh(2\eta J_3)}{\sinh(\eta)},   \label{eq:sl(2)-def}
\end{equation}
where $[x]_q \equiv \left(q^x - q^{-x}\right)/(q - q^{-1})$.\par
%
%
Such relations are left invariant under the transformation $\sigma\left(J_3\right)
= - J_3$, $\sigma\left(J_+\right) = J_-$, $\sigma\left(J_-\right) = J_+$. Hence,
defining ${\cal C} = \{+1, -1\}$, we get a finite, abelian colour group ${\cal G} =
\{\sigma^+ = \id, \sigma^- = \sigma\}$, isomorphic to the symmetric group~$S_2$.
The action of $\sigma^{\nu}$, $\nu = \pm1$, on the generators can be written in
compact form as 
\begin{equation}
  \sigma^{\nu}\left(J_3\right) = \nu J_3, \qquad \sigma^{\nu}\left(J_{\pm}\right) =
  J_{\pm\nu}.    \label{eq:S_2}
\end{equation}
\par
%
%
By using the results of Refs.~\cite{drinfeld, kirillov}, and
Definitions~\ref{def-colmaps} and~\ref{def-colR} of the present paper, we obtain
the following coloured comultiplication, counit, antipode, and universal $\cal
R$-matrix:
\begin{eqnarray}
  \Delta^{\lambda,\mu}_{q,\nu}\left(J_3\right) & = & (\lambda \nu)\, J_3 \otimes 1
           + (\mu \nu)\, 1 \otimes J_3, \nonumber \\
  \Delta^{\lambda,\mu}_{q,\nu}\left(J_{\pm}\right) & = & J_{\pm\lambda\nu}
           \otimes q^{\mu J_3} + q^{-\lambda J_3} \otimes J_{\pm\mu\nu}, \nonumber
           \\
  \epsilon_{q,\nu}(X) & = & 0, \qquad X \in \{J_3, J_{\pm}\}, \nonumber \\
  S^{\mu}_{q,\nu}\left(J_3\right) & = & - \mu \nu J_3, \qquad
           S^{\mu}_{q,\nu}\left(J_{\pm}\right) = - q^{\pm\nu} J_{\pm\mu\nu},
           \nonumber \\
  {\cal R}^{\lambda,\mu}_q & = & q^{2\lambda\mu J_3 \otimes J_3}
           \sum_{n=0}^{\infty} \frac{\left(1 - q^{-2}\right)^n}{[n]_q!}\, q^{n(n-1)/2}
           \left(q^{\lambda J_3} J_{\lambda}\right)^n \otimes \left(q^{-\mu J_3}
           J_{-\mu}\right)^n,
\end{eqnarray}
where $[n]_q! \equiv [n]_q [n-1]_q \ldots [1]_q$ for $n \in \N^+$, and $[0]_q! \equiv
1$.\par
%
%
The matrix representation of the coloured universal $\cal R$-matrix in any
finite-dimensional representation of $U_q(sl(2))$ provides us with a matrix
solution $R^{\lambda,\mu}_q$ of the coloured YBE~(\ref{eq:colYBE}), corresponding
to discrete colour parameters $\lambda$, $\mu = \pm1$. For instance, in the
two-dimensional representation of $U_q(sl(2))$,
\begin{equation}
  D(J_3) = \case{1}{2} \left(\begin{array}{cc}
                1 & 0 \\ 0 & -1
                \end{array} \right), \qquad 
  D(J_+) = \left(\begin{array}{cc}
                0 & 1 \\ 0 & 0
                \end{array} \right), \qquad
  D(J_-) = \left(\begin{array}{cc}
                0 & 0 \\ 1 & 0
                \end{array} \right),  \label{eq:sl(2)-rep}
\end{equation}
we get a (renormalized) $4 \times 4$ coloured $R$-matrix $R^{\lambda,\mu}_{q}
\equiv q^{1/2} (D \otimes D) \left({\cal R}^{\lambda,\mu}_q\right)$, whose
components are given by
\begin{eqnarray}
  R^{+,+}_q & = & \left(R^{-,-}_q\right)^t = R_q = \left(\begin{array}{cccc}
              q & 0 & 0              & 0 \\[0.1cm]
              0 & 1 & q - q^{-1} & 0 \\[0.1cm]
              0 & 0 & 1             & 0 \\[0.1cm]
              0 & 0 & 0             & q
              \end{array} \right), \nonumber \\
  R^{+,-}_q & = & \left(R^{-,+}_q\right)^t = \left(\begin{array}{cccc}
              1 & 0 & 0 & q - q^{-1} \\[0.1cm]
              0 & q & 0 & 0 \\[0.1cm]
              0 & 0 & q & 0 \\[0.1cm]
              0 & 0 & 0 & 1
              \end{array} \right),   \label{eq:sl(2)-R1} 
\end{eqnarray}
where $t$ stands for matrix transposition.\par
%
%
\subsubsection{The colour group ${\cal G} = Gl(1,\C)$}  
\label{subsubsec:Gl(1,C)}
The commutation relations~(\ref{eq:sl(2)-def}) are also left invariant under the
transformations
\begin{equation}
  \sigma^{\nu}\left(J_3\right) = J_3, \qquad  \sigma^{\nu}\left(J_{\pm}\right) =
  \nu^{\pm1} J_{\pm},   \label{eq:Gl(1,C)}
\end{equation}
where $\nu \in {\cal C} = \C \! \setminus \! \{0\}$. Since $\nu' \circ \nu = \nu' \nu$,
$\nu^0 = 1$, $\nu^i = \nu^{-1}$, the colour group~$\cal G$ is now isomorphic to the
abelian Lie group $Gl(1,\C)$.\par
%
%
The corresponding coloured quasitriangular Hopf algebra is defined
by~(\ref{eq:sl(2)-def}) and
\begin{eqnarray}
  \Delta^{\lambda,\mu}_{q,\nu}\left(J_3\right) & = & J_3 \otimes 1 + 1 \otimes J_3,
           \qquad \Delta^{\lambda,\mu}_{q,\nu}\left(J_{\pm}\right) =
           \left(\frac{\lambda}{\nu}\right)^{\pm1} J_{\pm} \otimes q^{J_3} +
           \left(\frac{\mu}{\nu}\right)^{\pm1} q^{- J_3} \otimes J_{\pm},
           \nonumber \\
  \epsilon_{q,\nu}(X) & = & 0, \qquad X \in \{J_3, J_{\pm}\}, \nonumber \\
  S^{\mu}_{q,\nu}\left(J_3\right) & = & - J_3, \qquad
           S^{\mu}_{q,\nu}\left(J_{\pm}\right) = - \left(\frac{\mu q}{\nu}
           \right)^{\pm1} J_{\pm}, \nonumber \\
  {\cal R}^{\lambda,\mu}_q & = & q^{2 J_3 \otimes J_3} \sum_{n=0}^{\infty}
           \frac{\left(1 - q^{-2}\right)^n}{[n]_q!}\, q^{n(n-1)/2}
           \left(\lambda q^{J_3} J_+\right)^n \otimes \left(\mu^{-1} q^{- J_3} J_-
           \right)^n.    \label{eq:sl(2)-Hopf}
\end{eqnarray}
\par
%
%
In the two-dimensional representation~(\ref{eq:sl(2)-rep}) of $U_q(sl(2))$, we
obtain a $4 \times 4$ matrix solution of the coloured YBE,
\begin{equation}
  R^{\lambda,\mu}_{q} \equiv q^{1/2} (D \otimes D) \left({\cal
  R}^{\lambda,\mu}_q\right) = \left(\begin{array}{cccc}
              q & 0 & 0                                                              & 0 \\[0.1cm]
              0 & 1 & \lambda \mu^{-1} \left(q - q^{-1}\right) & 0 \\[0.1cm]
              0 & 0 & 1                                                             & 0 \\[0.1cm]
              0 & 0 & 0                                                             & q
              \end{array} \right),   \label{eq:sl(2)-R2}
\end{equation}
depending upon continuous colour parameters $\lambda$, $\mu \in \C \! \setminus \!
\{0\}$. Similar results can be derived for other finite-dimensional
representations.\par
%
%
\subsubsection{The colour group ${\cal G} = Gl(1,\C) \sdprod S_2$}  
\label{subsubsec:sdprod}
The automorphisms, considered in Subsubsecs.~\ref{subsubsec:S_2}
and~\ref{subsubsec:Gl(1,C)}, can be combined by defining $\sigma^{\nu} \equiv
\sigma^{(\nu_1,\nu_2)} \equiv \sigma^{(\nu_1,+)} \circ \sigma^{(1,\nu_2)}$, where
$\sigma^{(\nu_1,+)}$ and $\sigma^{(1,\nu_2)}$ are given by Eqs.~(\ref{eq:Gl(1,C)})
and (\ref{eq:S_2}), respectively. Hence the colour set is the cartesian product ${\cal
C} = \left(\C \! \setminus \! \{0\}\right) \times \{+1,-1\}$,
\begin{equation}
  \sigma^{\nu}\left(J_3\right) = \nu_2 J_3, \qquad
  \sigma^{\nu}\left(J_{\pm}\right) = \nu_1^{\pm\nu_2} J_{\pm\nu_2},
\end{equation}
and $\nu' \circ \nu \equiv (\nu'_1, \nu'_2) \circ (\nu_1, \nu_2) = \left(\nu'_1
\nu_1^{\nu'_2}, \nu'_2 \nu_2\right)$, $\nu^0 = (1,+)$, $\nu^i = \left(\nu_1^{-\nu_2},
\nu_2\right)$. The colour group~$\cal G$ is nonabelian, and is a semidirect product
group, $Gl(1,\C) \sdprod S_2$. The subgroup $Gl(1,\C)$ is indeed invariant, whereas
$S_2$ is not, since $\nu' \circ (\nu_1,+) \circ \nu^{\prime\, i} =
\left(\nu_1^{\nu'_2},+\right)$, but $\nu' \circ (1,\nu_2) \circ \nu^{\prime\, i} =
\left((\nu'_1)^{1-\nu_2}, \nu_2\right)$.
\par
%
%
Eqs.~(\ref{eq:sl(2)-Hopf}) and~(\ref{eq:sl(2)-R2}) are now replaced by
\begin{eqnarray}
  \Delta^{\lambda,\mu}_{q,\nu}\left(J_3\right) & = & (\lambda_2 \nu_2)\, J_3
           \otimes 1 + (\mu_2 \nu_2)\, 1 \otimes J_3, \nonumber \\
  \Delta^{\lambda,\mu}_{q,\nu}\left(J_{\pm}\right) & = & \lambda_1^{\pm
           \lambda_2 \nu_2} \nu_1^{\mp1} J_{\pm\lambda_2\nu_2} \otimes q^{\mu_2
           J_3} + \mu_1^{\pm\mu_2\nu_2} \nu_1^{\mp1} q^{-\lambda_2 J_3} \otimes
           J_{\pm\mu_2\nu_2}, \nonumber \\
  \epsilon_{q,\nu}(X) & = & 0, \qquad X \in \{J_3, J_{\pm}\}, \nonumber \\
  S^{\mu}_{q,\nu}\left(J_3\right) & = & - \mu_2 \nu_2 J_3, \qquad
           S^{\mu}_{q,\nu}\left(J_{\pm}\right) = - \mu_1^{\pm\mu_2\nu_2}
           \nu_1^{\mp1} q^{\pm\nu_2} J_{\pm\mu_2\nu_2}, \nonumber \\
  {\cal R}^{\lambda,\mu}_q & = & q^{2\lambda_2\mu_2 J_3 \otimes J_3}
           \sum_{n=0}^{\infty} \frac{\left(1 - q^{-2}\right)^n}{[n]_q!}\, q^{n(n-1)/2}
           \left(\lambda_1^{\lambda_2} q^{\lambda_2 J_3} J_{\lambda_2}\right)^n
           \nonumber \\
  & & \mbox{} \otimes \left(\mu_1^{-\mu_2} q^{-\mu_2 J_3} J_{-\mu_2}\right)^n,
\end{eqnarray}
and
\begin{eqnarray}
  R^{(\lambda_1,+),(\mu_1,+)}_q & = &
              \left(R^{\left(\lambda_1^{-1},-\right),
              \left(\mu_1^{-1},-\right)}_q\right)^t = \left(\begin{array}{cccc}
              q & 0 & 0                                                                      & 0 \\[0.1cm]
              0 & 1 & \lambda_1 \mu_1^{-1} \left(q - q^{-1}\right) & 0 \\[0.1cm]
              0 & 0 & 1                                                                      & 0 \\[0.1cm]
              0 & 0 & 0                                                                      & q
              \end{array} \right), \nonumber \\
  R^{(\lambda_1,+),(\mu_1,-)}_q & = &
              \left(R^{\left(\lambda_1^{-1},-\right),
              \left(\mu_1^{-1},+\right)}_q\right)^t = \left(\begin{array}{cccc}
              1 & 0 & 0 & \lambda_1 \mu_1 \left(q - q^{-1}\right) \\[0.1cm]
              0 & q & 0 & 0 \\[0.1cm]
              0 & 0 & q & 0 \\[0.1cm]
              0 & 0 & 0 & 1
              \end{array} \right),   \label{eq:sl(2)-R3} 
\end{eqnarray}
respectively. The $4 \times 4$ matrices defined in Eqs.~(\ref{eq:sl(2)-R1}),
(\ref{eq:sl(2)-R2}), and~(\ref{eq:sl(2)-R3}) are known five-vertex solutions of
Eq.~(\ref{eq:colYBE})~\cite{hlavaty}.\par
%
%
\subsection{The two-parameter quantum algebra $U_{q,s}(gl(2))$} 
\label{subsec:gl(2)}
The next example deals with the two-parameter deformation of
$U(gl(2))$~\cite{schirrmacher}, whose universal $\cal R$-matrix was given in
Ref.~\cite{burdik92b}. Such an example is quite significant since $U_{q,s}(gl(2))$
and its corresponding quantum group have played an important role both in
generating some matrix solutions of the coloured YBE~\cite{ burdik92a,basu94},
and in constructing a coloured extension of the FRT formalism~\cite{basu94,
kundu94, basu95}.\par
%
%
The quantum algebra $U_{q,s}(gl(2))$, for which $k = \C$ and $q$, $s \in \C
\! \setminus \! \{0\}$, is generated by four operators $J_3$, $J_{\pm}$, $Z$, with
commutation relations
\begin{equation}
  \left[J_3, J_{\pm}\right] = \pm J_{\pm}, \qquad \left[J_+, J_-\right] = 
  \left[2J_3\right]_q, \qquad \left[Z, J_3\right] = \left[Z, J_{\pm}\right] = 0,
  \label{eq:gl(2)-def}
\end{equation}
and coalgebra and antipode depending upon both parameters $q$ and~$s$ (hence
${\cal Q} = \{\,(q,s)\,\}$).\par
%
%
Eq.~(\ref{eq:gl(2)-def}) is left invariant under the transformations
\begin{equation}
  \sigma^{\nu}\left(J_3\right) = J_3, \qquad  \sigma^{\nu}\left(J_{\pm}\right) =
  J_{\pm}, \qquad  \sigma^{\nu}\left(Z\right) = \nu Z,
\end{equation}
where $\nu \in {\cal C} = \C \! \setminus \! \{0\}$. As in
Subsubsec.~\ref{subsec:sl(2)}.\ref{subsubsec:Gl(1,C)}, the colour group is therefore
${\cal G} = Gl(1,\C)$. The coloured maps and universal $\cal R$-matrix are easily
obtained as
\begin{eqnarray}
  \Delta^{\lambda,\mu}_{q,s,\nu}\left(J_3\right) & = & J_3 \otimes 1 + 1 \otimes
           J_3, \qquad \Delta^{\lambda,\mu}_{q,s,\nu}\left(Z\right) = \frac{\lambda}
           {\nu}\, Z \otimes 1 + \frac{\mu}{\nu}\, 1 \otimes Z, \nonumber \\
  \Delta^{\lambda,\mu}_{q,s,\nu}\left(J_{\pm}\right) & = & J_{\pm} \otimes q^{J_3}
           \left(\frac{s}{q}\right)^{\pm\mu Z} + q^{- J_3} (q s)^{\pm\lambda Z} \otimes
           J_{\pm}, \nonumber \\
  \epsilon_{q,s,\nu}(X) & = & 0, \qquad X \in \{J_3, J_{\pm}, Z\}, \nonumber \\
  S^{\mu}_{q,s,\nu}\left(J_3\right) & = & - J_3, \qquad
           S^{\mu}_{q,s,\nu}\left(Z\right) = - \frac{\mu}{\nu} Z, \qquad
           S^{\mu}_{q,s,\nu}\left(J_{\pm}\right) = - q^{\pm1} s^{\mp 2\mu Z} J_{\pm},
           \nonumber \\
  {\cal R}^{\lambda,\mu}_{q,s} & = & q^{2 \left(J_3 \otimes J_3 - \lambda Z
           \otimes J_3 + \mu J_3 \otimes Z\right)} \sum_{n=0}^{\infty}
           \frac{\left(1 - q^{-2}\right)^n}{[n]_q!}\, q^{n(n-1)/2}
           \left(q^{J_3} (q s)^{- \lambda Z} J_+\right)^n \nonumber \\
  & & \mbox{} \otimes \left(q^{- J_3} \left(\frac{s}{q}\right)^{\mu Z} J_-\right)^n.   
\end{eqnarray}
\par
%
%
By considering for instance the matrix elements of~${\cal R}^{\lambda,\mu}_{q,s}$
in the $2 \times 2$ defining representation of $U_{q,s}(gl(2))$, given by
Eq.~(\ref{eq:sl(2)-rep}) and
\begin{equation}
  D(Z) = \left(\begin{array}{cc}
                1 & 0 \\ 0 & 1
                \end{array} \right), 
\end{equation}
we obtain the following $4 \times 4$ matrix solution of the coloured YBE,
\begin{equation}
  R^{\lambda,\mu}_{q,s} \equiv q^{1/2} (D \otimes D) \left({\cal
              R}^{\lambda,\mu}_{q,s}\right) = \left(\begin{array}{cccc}
              q^{1-\lambda+\mu} & 0 & 0 & 0 \\[0.1cm]
              0 & q^{\lambda+\mu} & \left(q - q^{-1}\right) s^{-\lambda+\mu} & 0 
                   \\[0.1cm]
              0 & 0 & q^{-\lambda-\mu} & 0 \\[0.1cm]
              0 & 0 & 0 & q^{1+\lambda-\mu}
              \end{array} \right).   
\end{equation}
The latter coincides with the coloured $R$-matrix previously derived by
Burd\'\i k and Hellin\-ger~\cite{burdik92a} by considering $2 \times 2$
representations of $U_{q,s}(gl(2))$ characterized by different eigenvalues
$\lambda$,~$\mu$ of~$Z$.\par
%
%
\subsection{The three-parameter quantum algebra $U_{q,s_1,s_2}(sl(3) \oplus
u(1) \oplus u(1))$} \label{sec:sl(3)}
Similar considerations to those in Subsec.~\ref{subsec:gl(2)} can be carried
through for some other multiparametric QUEA's of nonsemisimple Lie algebras. As
an example, let us consider the three-parameter deformation of $U(sl(3) \oplus
u(1) \oplus u(1))$, constructed by Burd\'\i k and Hellinger~\cite{burdik92a}. It is
generated by eight operators $H_i$, $X^{\pm}_i$, $Z_i$, $i=1$,~2, with relations
\begin{eqnarray}
  \left[H_i, X^{\pm}_i\right] & = & \pm 2 X^{\pm}_i, \qquad \left[H_i,
           X^{\pm}_j\right] = \mp X^{\pm}_j \qquad (i\ne j), \nonumber \\
  \left[X^+_i, X^-_j\right] & = & \delta_{ij} \left[H_i\right]_q, \qquad
           \left[X^{\pm}_1, X^{\pm}_2\right] = 0, \nonumber \\
  q^{-1/2} X^{\pm}_1 X^{\pm}_3 - q^{1/2} X^{\pm}_3 X^{\pm}_1 & = & q^{1/2}
           X^{\pm}_2 X^{\pm}_3 - q^{-1/2} X^{\pm}_3 X^{\pm}_2 = 0, \nonumber \\  
  \left[Z_i, H_j\right] & = & \left[Z_i, X^{\pm}_j\right] = 0,  \label{eq:sl(3)-def}
\end{eqnarray}
where $X^{\pm}_3 \equiv q^{1/2} X^{\pm}_1 X^{\pm}_2 - q^{-1/2} X^{\pm}_2
X^{\pm}_1$. Here $k = \C$ and ${\cal Q} = \{\,(q, s_1, s_2)\,\}$, where $q$, $s_1$,
$s_2 \in \C \! \setminus \! \{0\}$, and $s_1$, $s_2$ make their appearance only in 
the coalgebra structure and the antipode.\par
%
%
Eq.~(\ref{eq:sl(3)-def}) is left invariant under the transformations
\begin{equation}
  \sigma^{\nu}\left(H_i\right) = H_i, \qquad  \sigma^{\nu}\left(X^{\pm}_i\right) =
  X^{\pm}_i, \qquad  \sigma^{\nu}\left(Z_i\right) = \nu_i Z_i,
\end{equation}
where $\nu \equiv (\nu_1, \nu_2) \in {\cal C} = (\C \! \setminus \! \{0\}) \times (\C
\! \setminus \! \{0\})$. Hence the colour group~$\cal G$ is a direct product group,
$Gl(1,\C) \otimes Gl(1,\C)$.\par
%
%
The coloured maps and universal $\cal R$-matrix are given by
\begin{eqnarray}
  \Delta^{\lambda,\mu}_{q,s_1,s_2,\nu}\left(H_i\right) & = & H_i \otimes 1 + 1
           \otimes H_i, \qquad \Delta^{\lambda,\mu}_{q,s_1,s_2,\nu}\left(Z_i\right) =
           \frac{\lambda_i}{\nu_i}\, Z_i \otimes 1 + \frac{\mu_i}{\nu_i}\, 1 \otimes
           Z_i, \nonumber \\
  \Delta^{\lambda,\mu}_{q,s_1,s_2,\nu}\left(X^{\pm}_i\right) & = & X^{\pm}_i
           \otimes q^{H_i/2} (q s_i)^{\pm\mu_i Z_i/2} + q^{-H_i/2}
           \left(\frac{s_i}{q}\right)^{\pm\lambda_i Z_i/2} \otimes X^{\pm}_i,
           \nonumber \\
  \epsilon_{q,s_1,s_2,\nu}(X) & = & 0, \qquad X \in \{H_i, X^{\pm}_i, Z_i\}, \nonumber
           \\
  S^{\mu}_{q,s_1,s_2,\nu}\left(H_i\right) & = & - H_i, \qquad
           S^{\mu}_{q,s_1,s_2,\nu}\left(Z_i\right) = - \frac{\mu_i}{\nu_i} Z_i, 
           \qquad S^{\mu}_{q,s_1,s_2,\nu}\left(X^{\pm}_i\right) = - q^{\pm1} 
           s_i^{\mp\mu_i Z_i} X^{\pm}_i, \nonumber \\
  {\cal R}^{\lambda,\mu}_{q,s_1,s_2} & = & q^{\sum_{i,j} a^{-1}_{i,j} \left(
           \lambda_i Z_i \otimes H_j - \mu_j H_i \otimes Z_j + H_i \otimes
           H_j\right)} E_{q^{-2}}\left(\alpha e^{\lambda}_1 \otimes f^{\mu}_1\right)
           \nonumber \\
  & & \mbox{} \times E_{q^{-2}}\left(-\alpha e^{\lambda}_3 \otimes
           f^{\mu}_3\right) E_{q^{-2}}\left(\alpha e^{\lambda}_2 \otimes
           f^{\mu}_2\right),   \label{eq:sl(3)-Hopf}  
\end{eqnarray}
where $\alpha \equiv 1 - q^{-2}$, $E_x(A)$ is the $q$-exponential
\begin{equation}
  E_x(A) = \sum_{n=0}^{\infty} \frac{x^{-n(n-1)/4}}{[n]_{x^{1/2}}!} A^n,
\end{equation}
$a$ is the Cartan matrix of $sl(3)$, $a^{-1}$ its inverse,
\begin{equation}
  a = \left(\begin{array}{cc}
        2   & -1 \\
        -1 & 2
       \end{array}\right), \qquad
  a^{-1} = \case{1}{3} \left(\begin{array}{cc}
        2   & 1 \\
        1   & 2
       \end{array}\right),
\end{equation}
and
\begin{eqnarray}
  e^{\lambda}_i & = & q^{H_i/2} \left(\frac{s_i}{q}\right)^{-\lambda_i Z_i/2} X^+_i,
  \qquad f^{\mu}_i = q^{-H_i/2} (q s_i)^{\mu_i Z_i/2} X^-_i, \nonumber \\
  e^{\lambda}_3 & = & e^{\lambda}_1 e^{\lambda}_2 - q^{-1} e^{\lambda}_2
  e^{\lambda}_1, \qquad f^{\mu}_3 = f^{\mu}_1 f^{\mu}_2 - q^{-1} f^{\mu}_2
  f^{\mu}_1.   
\end{eqnarray}
\par
%
%
In the $3 \times 3$ defining representation of $U_{q,s_1,s_2}(sl(3) \oplus u(1)
\oplus u(1))$, the coloured universal $\cal R$-matrix, given in
Eq.~(\ref{eq:sl(3)-Hopf}), gives rise to a $9 \times 9$ matrix solution of the
coloured YBE. As it coincides with the matrix given in Eqs.~(28) and~(29) of
Ref.~\cite{burdik92a}, we shall not reproduce it here.\par
%
%
\section{EXAMPLES OF COLOURED QUEA'S WITH VARYING PARAMETERS}
\label{sec:varying}
\setcounter{equation}{0}
The examples considered in the present section differ from those constructed in
Sec.~\ref{sec:fixed} by the fact that the parameter set~$\cal Q$ now contains
more than one element and the transformations of the colour group~$\cal G$ in
general change the parameters (hence $q^{\nu} \ne q$). Although we simultaneously
deal with elements~$a$ belonging to different Hopf algebras~${\cal H}_{q^{\nu}}$
of the set~$\cal H$, we denote them in the same way in order not to overload
notation by adding an extra index referring to the corresponding algebra, or colour
parameter. This should cause no confusion since, from the context, it is always
clear to which Hopf algebra $a$ belongs.\par
%
%
\subsection{The nonstandard quantum algebra $U_h(sl(2))$}  \label{sec:ns-sl(2)}
We begin by considering the nonstandard (Jordanian) deformation of
$U(sl(2))$~\cite{ohn}, known as $U_h(sl(2))$, $h \in \C \! \setminus \! \{0\}$, which 
can be obtained by contracting the standard DJ deformation
$U_q(sl(2))$~\cite{agha}. Its universal $\cal R$-matrix was derived in
Refs.~\cite{shariati, ballesteros96a}.\par
%
%
$U_h(sl(2))$ is generated by three operators $H$, $J_{\pm}$, satisfying the
commutation relations~\cite{ohn}
\begin{eqnarray}
  \left[H, J_+\right] & = & 2\, \frac{\sinh(h J_+)}{h}, \qquad \left[J_+, J_-\right] =
              H, \nonumber \\
  \left[H, J_-\right] & = & - J_- \cosh(h J_+) - \cosh(h J_+) J_-.  
              \label{eq:ns-sl(2)-def} 
\end{eqnarray}
Its universal $\cal R$-matrix assumes a very simple form~\cite{ballesteros96a}
provided one chooses another basis, whose generators $A$,~$A_{\pm}$ are defined
in terms of the old ones by
\begin{equation}
  A = e^{h J_+} J_3, \qquad A_+ = J_+, \qquad A_- = e^{h J_+} \left(J_- - \frac{h}{4}
  \sinh(h J_+)\right).
\end{equation}
Eq.~(\ref{eq:ns-sl(2)-def}) is then transformed into
\begin{equation}
  \left[A, A_+\right] = \frac{e^{2h A_+}-1}{h}, \qquad \left[A, A_-\right] = - 2A_-
  + h A^2, \qquad \left[A_+, A_-\right] = A.
\end{equation}
\par
%
%
Such relations are left invariant under the transformations
\begin{equation}
  \sigma^{\nu}(A) = A, \qquad \sigma^{\nu}\left(A_+\right) = \nu A_+, \qquad
  \sigma^{\nu}\left(A_-\right) = \nu^{-1} A_-,
\end{equation}
where $\nu \in {\cal C} = \C \! \setminus \! \{0\}$, with the proviso that $h$ 
becomes $h^{\nu} = \nu h$ (hence ${\cal Q} = \C \! \setminus \! \{0\}$). The
$\sigma^{\nu}$'s are therefore isomorphic mappings from $U_h(sl(2))$ to
$U_{h^{\nu}}(sl(2))$, and define a colour group ${\cal G} = Gl(1,\C)$. The
corresponding coloured maps and universal
$\cal R$-matrix are given by the relations
\begin{eqnarray}
  \Delta^{\lambda,\mu}_{h,\nu}\left(A_+\right) & = & \frac{\lambda}{\nu}\, A_+
            \otimes 1 + \frac{\mu}{\nu}\, 1 \otimes A_+, \qquad  
            \Delta^{\lambda,\mu}_{h,\nu}\left(A\right) = A \otimes e^{2\mu h A_+} 
            + 1 \otimes A, \nonumber \\
  \Delta^{\lambda,\mu}_{h,\nu}\left(A_-\right) & = & \frac{\nu}{\lambda}\, A_-
            \otimes e^{2\mu h A_+} + \frac{\nu}{\mu}\, 1 \otimes A_-, \nonumber \\
  \epsilon_{h,\nu}(X) & = & 0, \qquad x \in \{A, A_{\pm}\}, \nonumber \\
  S^{\mu}_{h,\nu}\left(A_+\right) & = & - \frac{\mu}{\nu} A_+, \qquad
            S^{\mu}_{h,\nu}\left(A\right) = - A e^{-2\mu h A_+}, \qquad 
            S^{\mu}_{h,\nu}\left(A_-\right) = - \frac{\nu}{\mu} A_- e^{-2\mu h A_+},
            \nonumber \\
  {\cal R}^{\lambda,\mu}_h & = & \exp\left\{- \lambda h A_+ \otimes A\right\}
            \exp\left\{\mu h A \otimes A_+\right\}.  \label{eq:ns-sl(2)-Hopf}  
\end{eqnarray}
\par
%
%
The two-dimensional defining representation of $U_h(sl(2))$ is still given by
Eq.~(\ref{eq:sl(2)-rep}), where $D(H) = 2 D(J_3)$, thence~\cite{footnote2}
\begin{equation}
  D(A) = \left(\begin{array}{cc}
                1 & -h \\ 0 & -1
                \end{array} \right), \qquad 
  D(A_+) = \left(\begin{array}{cc}
                0 & 1 \\ 0 & 0
                \end{array} \right), \qquad
  D(A_-) = \left(\begin{array}{cc}
                h & - h^2/4 \\ 1 & 0
                \end{array} \right).   \label{eq:ns-sl(2)-rep}
\end{equation}
From Eq.~(\ref{eq:ns-sl(2)-Hopf}), it is straightforward to obtain
\begin{equation}
  R^{\lambda,\mu}_h \equiv (D \otimes D)\left({\cal R}^{\lambda,\mu}_h\right) =
  \left(\begin{array}{cccc}
        1 & \mu h & - \lambda h & (\lambda - \mu + \lambda \mu) h^2 \\[0.1cm]
        0 & 1        & 0                 & \lambda h \\[0.1cm]
        0 & 0        & 1                 & - \mu h \\[0.1cm]
        0 & 0        & 0                 & 1    
  \end{array}\right).
\end{equation}
To the best of our knowledge, this $4 \times 4$ matrix solution of the coloured YBE
is new. Higher-dimensional solutions of the latter could be obtained in a similar
way by considering other finite-dimensional irreducible representations of the
quantum algebra $U_h(sl(2))$~\cite{ballesteros96a, abdesselam}.\par
%
%
\subsection{The standard quantum oscillator algebra $U^{(s)}_z(h(4))$}
\label{subsec:s-h(4)}
The next example is the standard deformation $U^{(s)}_z(h(4))$ of the oscillator
algebra $U(h(4))$, which was first derived by contracting
$U_q(gl(2))$~\cite{celeghini90}, then recently obtained in a more convenient
basis~\cite{ballesteros96b} by using the Lyakhovsky and Mudrov
formalism~\cite{lyakhovsky}. This quantum algebra has been used to construct
a solution of the coloured YBE connected with some link
invariants~\cite{gomez}.\par
%
%
$U^{(s)}_z(h(4))$ is generated by four operators $N$, $M$,~$A_{\pm}$, satisfying
the commutation relations
\begin{equation}
  \left[N, A_{\pm}\right] = \pm A_{\pm}, \qquad \left[A_-, A_+\right] =
  \frac{\sinh(zM)}{z}, \qquad \left[M, N\right] = \left[M, A_{\pm}\right] = 0.
  \label{eq:s-h(4)-def}
\end{equation}
Here we assume $k = \C$ and $z \in {\cal Q} = \C \! \setminus \! \{0\}$. The algebra
defining relations~(\ref{eq:s-h(4)-def}) are left invariant under the
transformations
\begin{equation}
  \sigma^{\nu}(N) = N, \qquad \sigma^{\nu}(M) = \nu_+ \nu_- M, \qquad
  \sigma^{\nu}(A_+) = \nu_+ A_+, \qquad \sigma^{\nu}(A_-) = \nu_- A_-,
  \label{eq:s-h(4)-G}
\end{equation}
where $\nu \equiv (\nu_+, \nu_-)$, provided $z$ is changed into $z^{\nu} = \nu_+
\nu_- z$. Hence the colour set is the cartesian product ${\cal C} = (\C \! \setminus \!
\{0\}) \times (\C \! \setminus \! \{0\})$, and the colour group is the direct product
group ${\cal G} = Gl(1,\C) \otimes Gl(1,\C)$.\par
%
%
The corresponding coloured maps and universal $\cal R$-matrix are given by
\begin{eqnarray}
  \Delta^{\lambda,\mu}_{z,\nu}(N) & = & N \otimes 1 + 1 \otimes N, \qquad
          \Delta^{\lambda,\mu}_{z,\nu}(M) = \frac{\lambda_+\lambda_-}{\nu_+\nu_-}\,
          M \otimes 1 + \frac{\mu_+\mu_-}{\nu_+\nu_-}\, 1 \otimes M, \nonumber \\
  \Delta^{\lambda,\mu}_{z,\nu}(A_+) & = & \frac{\lambda_+}{\nu_+}\, A_+ \otimes 1 
          + \frac{\mu_+}{\nu_+}\, e^{-\lambda_+\lambda_- z M} \otimes A_+,
          \nonumber \\
  \Delta^{\lambda,\mu}_{z,\nu}(A_-) & = & \frac{\lambda_-}{\nu_-}\, A_- \otimes  
          e^{\mu_+\mu_- z M}+ \frac{\mu_-}{\nu_-}\, 1 \otimes A_-, \nonumber \\
  \epsilon_{z,\nu}(X) & = & 0, \qquad X \in \{N, M, A_{\pm}\}, \nonumber \\
  S^{\mu}_{z,\nu}(N) & = & - N, \qquad S^{\mu}_{z,\nu}(M) = -
          \frac{\mu_+\mu_-}{\nu_+\nu_-} M, \qquad S^{\mu}_{z,\nu}(A_{\pm}) = -
          \frac{\mu_{\pm}}{\nu_{\pm}} A_{\pm}\, e^{\pm\mu_+\mu_- z M}, \nonumber
          \\
  {\cal R}^{\lambda,\mu}_z & = & \exp\{-\lambda_+ \lambda_- z M \otimes N\}
          \exp\{-\mu_+ \mu_- z N \otimes M\} \nonumber \\
  & & \mbox{} \times \exp\{2\lambda_- \mu_+ z A_- \otimes A_+\}. 
          \label{eq:s-h(4)-Hopf}
\end{eqnarray}
\par
%
%
In the $3 \times 3$ matrix representation of $U^{(s)}_z(h(4))$ defined by
\begin{eqnarray}
  D(N) & = & \left(\begin{array}{ccc}
        0 & 0 & 0 \\ 0 & 1 & 0 \\ 0 & 0 & 0
        \end{array}\right), \qquad
  D(M) = \left(\begin{array}{ccc}
        0 & 0 & 1 \\ 0 & 0 & 0 \\ 0 & 0 & 0
        \end{array}\right), \nonumber \\
  D(A_+) & = & \left(\begin{array}{ccc}
        0 & 0 & 0 \\ 0 & 0 & 1 \\ 0 & 0 & 0
        \end{array}\right), \qquad
  D(A_-) = \left(\begin{array}{ccc}
        0 & 1 & 0 \\ 0 & 0 & 0 \\ 0 & 0 & 0
        \end{array}\right),   \label{eq:s-h(4)-rep}
\end{eqnarray}
the coloured universal $\cal R$-matrix is represented by the $9 \times 9$ matrix
\begin{equation}
  R^{\lambda,\mu}_z \equiv (D \otimes D)\left({\cal R}^{\lambda,\mu}_z\right) =
  \left(\begin{array}{ccc}
          1_3 & 2\lambda_-\mu_+ z D(A_+) & -\lambda_+\lambda_- z D(N) \\[0.1cm]
          0_3 & 1_3 - \mu_+\mu_- z D(M) & 0_3 \\[0.1cm]
          0_3 & 0_3 & 1_3
  \end{array}\right),   \label{s-h(4)-R}
\end{equation}
where $1_3$ and $0_3$ denote the $3\times 3$ unit and null matrices
respectively.\par
%
%
{\it Remarks.\/} (1) Since in physical applications, $A_+$ (resp.~$A_-$) is
interpreted as a creation (resp.~annihilation) operator, and $N$ as a number
operator, one actually deals there with a real form of $U^{(s)}_z(h(4))$,
corresponding to the star operation (Hermitian conjugation)
\begin{equation}
  N^{\dagger} = N, \qquad M^{\dagger} = M, \qquad A_{\pm}^{\dagger} = A_{\mp},
  \label{eq:s-h(4)-star}
\end{equation}
and to real or imaginary values of~$z$. Eq.~(\ref{eq:s-h(4)-star}) restricts
$\nu_{\pm}$ in Eq.~(\ref{eq:s-h(4)-G}) to values satisfying the condition $\nu_- =
\overline{\nu_+}$, where the bar denotes complex conjugation. Hence, in such a
case, ${\cal C} = \C \! \setminus \! \{0\}$, and ${\cal G} = Gl(1,\C)$. (2) It is easy to
endow $U^{(s)}_z(h(4))$ with a coloured Hopf structure corresponding to a
nonabelian colour group by combining transformations~(\ref{eq:s-h(4)-G}) with the
elements of~$S_2$, defined by $\sigma^+ = \id$, and $\sigma^-(N) = - N$,
$\sigma^-(M) = - M$, $\sigma^-(A_{\pm}) = A_{\mp}$ (implying $z^{\pm} = z$). This
can be done along the same lines as in the example discussed in
Subsubsec.~\ref{sec:fixed}.\ref{subsec:sl(2)}.\ref{subsubsec:sdprod}.\par
%
%
\subsection{The one-parameter nonstandard quantum oscillator algebra
$U^{(n)}_z(h(4))$}
\label{subsec:1ns-h(4)}
Instead of the standard deformation of the oscillator algebra, dealt with in
Subsec.~\ref{subsec:s-h(4)}, we consider here the one-parameter nonstandard type
I$_+$ deformation of the same, constructed in Ref.~\cite{ballesteros96b}, and
denoted there by $U^{(n)}_z(h(4))$. For such an algebra, Eq.~(\ref{eq:s-h(4)-def}) is
replaced by
\begin{eqnarray}
  \left[N, A_+\right] & = & \frac{e^{z A_+} - 1}{z}, \qquad \left[N, A_-\right] = - A_-,
             \qquad \left[A_-, A_+\right] = M e^{z A_+}, \nonumber \\
  \left[M, N\right] & = & \left[M, A_{\pm}\right] = 0,  \label{eq:1ns-h(4)-def}
\end{eqnarray}
where we choose $z \in {\cal Q} = \C \! \setminus \! \{0\}$.\par
%
%
The defining relations~(\ref{eq:1ns-h(4)-def}) are left invariant under the
transformations
\begin{equation}
  \sigma^{\nu}(N) = N, \qquad \sigma^{\nu}(M) = \nu M, \qquad
  \sigma^{\nu}(A_+) = \nu A_+, \qquad \sigma^{\nu}(A_-) = A_-, 
  \label{eq:1ns-h(4)-G}
\end{equation}
provided $z$ is changed into $z^{\nu} = \nu z$. Hence ${\cal C} = \C \! \setminus \!
\{0\}$, and ${\cal G} = Gl(1,\C)$.\par
%
%
The counterparts of Eqs.~(\ref{eq:s-h(4)-Hopf}) and~(\ref{s-h(4)-R}) are now
\begin{eqnarray}
  \Delta^{\lambda,\mu}_{z,\nu}(N) & = & N \otimes e^{\mu z A_+} + 1 \otimes N,
          \qquad \Delta^{\lambda,\mu}_{z,\nu}(M) = \frac{\lambda}{\nu}\, M \otimes 1 +
          \frac{\mu}{\nu}\, 1 \otimes M, \nonumber \\
  \Delta^{\lambda,\mu}_{z,\nu}(A_+) & = & \frac{\lambda}{\nu}\, A_+ \otimes 1 
          + \frac{\mu}{\nu}\, 1 \otimes A_+, \nonumber \\
  \Delta^{\lambda,\mu}_{z,\nu}(A_-) & = & A_- \otimes  e^{\mu z A_+} + 1 \otimes
          A_- + \mu z N \otimes M e^{\mu z A_+}, \nonumber \\
  \epsilon_{z,\nu}(X) & = & 0, \qquad X \in \{N, M, A_{\pm}\}, \nonumber \\
  S^{\mu}_{z,\nu}(N) & = & - N e^{-\mu z A_+}, \qquad S^{\mu}_{z,\nu}(M) = -
          \frac{\mu}{\nu} M, \qquad S^{\mu}_{z,\nu}(A_+) = - \frac{\mu}{\nu} A_+,
          \nonumber \\
  S^{\mu}_{z,\nu}(A_-) & = & -A_- e^{-\mu z A_+} + \mu z N M e^{-\mu z A_+},
          \nonumber \\
  {\cal R}^{\lambda,\mu}_z & = & \exp\{-\lambda z A_+ \otimes N\}
          \exp\{\mu z N \otimes A_+\},   \label{eq:1ns-h(4)-Hopf}
\end{eqnarray}
and
\begin{equation}
  R^{\lambda,\mu}_z \equiv (D \otimes D)\left({\cal R}^{\lambda,\mu}_z\right) =
  \left(\begin{array}{ccc}
          1_3 & 0_3 & 0_3 \\[0.1cm]
          0_3 & 1_3 + \mu z D(A_+) & -\lambda z D(N) \\[0.1cm]
          0_3 & 0_3 & 1_3
  \end{array}\right),  \label{eq:1ns-h(4)-R}
\end{equation}
respectively, where D is again defined by Eq.~(\ref{eq:s-h(4)-rep}).\par
%
%
{\it Remarks.\/} (1) Eq.~(\ref{eq:1ns-h(4)-def}) is not compatible with the star
operation~(\ref{eq:s-h(4)-star}). (2) In Eq.~(\ref{eq:1ns-h(4)-G}), we could multiply
both $M$ and~$A_-$ by some extra parameter, but this would not modify the
coloured universal $\cal R$-matrix, as given in Eq.~(\ref{eq:1ns-h(4)-Hopf}).\par
%
%
\subsection{The three-parameter nonstandard quantum oscillator algebra
$U^{(IIn)}_{\vartheta,\beta_+,\beta_-}(h(4))$}
\label{subsec:3ns-h(4)}
Similar results can be obtained for more complicated deformations of the oscillator
algebra. In the present subsection, we consider the three-parameter nonstandard
deformation constructed in Ref.~\cite{ballesteros96b}, where it is denoted by
$U^{(IIn)}_{\vartheta,\beta_+,\beta_-}(h(4))$. The algebra defining relations are 
\begin{eqnarray}
  \left[N, A_+\right] & = & A_+ - \beta_- V(-\vartheta), \qquad 
               \left[N, A_-\right] = - A_- - \beta_+ V(\vartheta), \nonumber \\
  \left[A_-, A_+\right] & = & M, \qquad \left[M, N\right] = \left[M, A_{\pm}\right] =
               0,    \label{eq:3ns-h(4)-def}
\end{eqnarray}
where
\begin{equation}
  V(x) \equiv \frac{1}{x^2} \left(e^{xM} - 1 - x M\right),
\end{equation}
and we assume $(\vartheta, \beta_+, \beta_-) \in {\cal Q} = (\C \! \setminus \! \{0\})
\times (\C \! \setminus \! \{0\}) \times (\C \! \setminus \! \{0\})$.\par
%
%
The transformations
\begin{equation}
  \sigma^{\nu}(N) = N, \qquad \sigma^{\nu}(M) = \nu_+ \nu_- M, \qquad
  \sigma^{\nu}(A_+) = \nu_+ A_+, \qquad \sigma^{\nu}(A_-) = \nu_- A_-,
\end{equation}
where $\nu \equiv (\nu_+, \nu_-)$, leave Eq.~(\ref{eq:3ns-h(4)-def}) invariant,
provided $\vartheta$, $\beta_+$,~$\beta_-$ are changed into $\vartheta^{\nu} =
\nu_+ \nu_- \vartheta$, $\beta_+^{\nu} = \nu_+^2 \nu_- \beta_+$, $\beta_-^{\nu} =
\nu_+ \nu_-^2 \beta_-$, respectively. Hence ${\cal C} = (\C \! \setminus \! \{0\}) 
\times (\C \! \setminus \! \{0\})$, and ${\cal G} = Gl(1,\C) \otimes Gl(1,\C)$.\par
%
%
The counterparts of Eqs.~(\ref{eq:s-h(4)-Hopf}) and~(\ref{s-h(4)-R}) are now
\begin{eqnarray}
  \Delta^{\lambda,\mu}_{\vartheta,\beta_+,\beta_-,\nu}(N) & = & N \otimes 1 + 1
          \otimes N + \frac{\lambda_+ \beta_+}{\vartheta}\, A_+ \otimes \left(1 -
          e^{-\mu_+ \mu_- \vartheta M}\right) \nonumber \\
  & & \mbox{} + \frac{\lambda_- \beta_-}{\vartheta}\, A_- \otimes \left(1 -
          e^{\mu_+ \mu_- \vartheta M}\right), \nonumber \\
  \Delta^{\lambda,\mu}_{\vartheta,\beta_+,\beta_-,\nu}(M) & = &
          \frac{\lambda_+\lambda_-}{\nu_+\nu_-}\, M \otimes 1 +
          \frac{\mu_+\mu_-}{\nu_+\nu_-}\, 1 \otimes M, \nonumber \\
  \Delta^{\lambda,\mu}_{\vartheta,\beta_+,\beta_-,\nu}(A_{\pm}) & = &
          \frac{\lambda_{\pm}}{\nu_{\pm}}\, A_{\pm} \otimes e^{\mp\mu_+\mu_-
          \vartheta M} + \frac{\mu_{\pm}}{\nu_{\pm}}\, 1 \otimes A_{\pm}, \nonumber
          \\
  \epsilon_{\vartheta,\beta_+,\beta_-,\nu}(X) & = & 0, \qquad X \in \{N, M,
         A_{\pm}\}, \nonumber \\
  S^{\mu}_{\vartheta,\beta_+,\beta_-,\nu}(N) & = & - N - \frac{\mu_+
         \beta_+}{\vartheta} A_+  \left(1 - e^{\mu_+ \mu_- \vartheta M}\right) 
         - \frac{\mu_- \beta_-}{\vartheta} A_-  \left(1 - e^{-\mu_+ \mu_- \vartheta
         M}\right), \nonumber \\
  S^{\mu}_{\vartheta,\beta_+,\beta_-,\nu}(M) & = & - \frac{\mu_+\mu_-}
         {\nu_+\nu_-} M, \qquad S^{\mu}_{\vartheta,\beta_+,\beta_-,\nu}(A_{\pm}) = -
         \frac{\mu_{\pm}}{\nu_{\pm}} A_{\pm} e^{\pm\mu_+\mu_- \vartheta M},
         \nonumber \\
  {\cal R}^{\lambda,\mu}_{\vartheta,\beta_+,\beta_-} & = & \exp\{-\lambda_+
         \lambda_- M \otimes (\vartheta N + \mu_+ \beta_+ A_+ + \mu_- \beta_-
         A_-)\} \nonumber \\
  & & \mbox{} \times \exp\{\mu_+ \mu_- (\vartheta N + \lambda_+ \beta_+ A_+ +
         \lambda_- \beta_- A_-) \otimes M\},
\end{eqnarray}
and
\begin{eqnarray}
  R^{\lambda,\mu}_{\vartheta,\beta_+,\beta_-} & \equiv & (D \otimes
           D)\left({\cal R}^{\lambda,\mu}_{\vartheta,\beta_+,\beta_-}\right)
            \nonumber \\
  & = &\left(\begin{array}{ccc}
          \ss 1_3 & \ss \lambda_-\mu_+\mu_-\beta_- D(M) & \ss
                -\lambda_+\lambda_- \left(\vartheta D(N) + \mu_+\beta_+ D(A_+) +
                \mu_-\beta_- D(A_-)\right) \\
          \ss 0_3 & \ss 1_3 + \mu_+\mu_-\vartheta D(M) & \ss
                \lambda_+\mu_+\mu_-\beta_+ D(M) \\
          \ss 0_3 & \ss 0_3 & \ss 1_3
  \end{array}\right),   \label{eq:3ns-h(4)-R}
\end{eqnarray}
respectively, where $D$ is again defined by Eq.~(\ref{eq:s-h(4)-rep}).
Eq.~(\ref{eq:3ns-h(4)-R}), as well as Eq.~(\ref{eq:1ns-h(4)-R}), provide new
matrix solutions of the coloured YBE.\par
%
%
{\it Remark.\/} The real form of $U^{(IIn)}_{\vartheta,\beta_+,\beta_-}(h(4))$,
corresponding to the star operation~(\ref{eq:s-h(4)-star}), is obtained for
$\vartheta = - \overline{\vartheta}$, and $\beta_- = - \overline{\beta_+}$. The
colour parameters are then restricted by the condition $\nu_- = \overline{\nu_+}$,
so that we are left with ${\cal C} = \C \! \setminus \! \{0\}$, and ${\cal G} =
Gl(1,\C)$.\par 
%
%
\subsection{The standard three-dimensional quantum Euclidean algebra
$U_w(e(3))$}    \label{subsec:euclide}
In the present subsection, we consider the three-dimensional quantum Euclidean 
algebra $U_w(e(3))$, which was obtained by contracting the standard DJ
deformation of $so(4)$~\cite{celeghini91}.\par
%
%
A basis of $U_w(e(3))$ is made of six operators $J_3$, $J_{\pm}$,
$P_3$,~$P_{\pm}$, generating rotations and translations in the $w \to 0$ limit
respectively, and satisfying the commutation relations
\begin{eqnarray}
  \left[J_3, J_{\pm}\right] & = & \pm J_{\pm}, \qquad  \left[J_+, J_-\right] = 2J_3
           \cosh(2wP_3), \nonumber \\
  \left[J_3, P_{\pm}\right] & = & \left[P_3, J_{\pm}\right] = \pm P_{\pm}, \qquad
           \left[J_{\pm}, P_{\mp}\right] = \pm \frac{\sinh(2wP_3)}{w}, \nonumber \\
  \left[J_3, P_3\right] & = & \left[J_{\pm}, P_{\pm}\right] = \left[P_3,
           P_{\pm}\right] = \left[P_+, P_-\right] = 0.  \label{eq:euclide-def}
\end{eqnarray} 
Here we assume $k = \R$, and $w \in {\cal Q} = \R \! \setminus \! \{0\}$, which is
compatible with the star operation usually imposed on $U(e(3))$, namely
\begin{equation}
  J_3^{\dagger} = J_3, \qquad  J_{\pm}^{\dagger} = J_{\mp}, \qquad  P_3^{\dagger} =
  P_3, \qquad  P_{\pm}^{\dagger} = P_{\mp}.  
\end{equation}
\par
%
%
The algebra defining relations~(\ref{eq:euclide-def}) are left invariant under the
transformations
\begin{equation}
  \sigma^{\nu}(J_3) = J_3, \qquad \sigma^{\nu}(J_{\pm}) = J_{\pm}, \qquad
  \sigma^{\nu}(P_3) = \nu P_3, \qquad \sigma^{\nu}(P_{\pm}) = \nu P_{\pm}, 
\end{equation}
where $\nu \in \R \! \setminus \! \{0\}$, provided $w$ is changed into $w^{\nu} = 
\nu w$. Hence the colour set and the colour group are ${\cal C} = \R \! \setminus \!
\{0\}$, and ${\cal G} = Gl(1,\R)$, respectively.\par
%
%
The coloured comultiplication, counit, antipode, and universal $\cal R$-matrix are
easily found to be given by
\begin{eqnarray}
  \Delta^{\lambda,\mu}_{w,\nu}(J_3) & = & J_3 \otimes 1 + 1 \otimes J_3, \qquad
           \Delta^{\lambda,\mu}_{w,\nu}(P_3) = \frac{\lambda}{\nu}\, P_3 \otimes 1
           + \frac{\mu}{\nu}\, 1 \otimes P_3, \nonumber \\
  \Delta^{\lambda,\mu}_{w,\nu}(J_{\pm}) & = & J_{\pm} \otimes e^{\mu w P_3} + 
            e^{-\lambda w P_3} \otimes J_{\pm} + w \left(\lambda P_{\pm} \otimes
            e^{\mu w P_3} J_3 - \mu e^{-\lambda w P_3} J_3 \otimes P_{\pm}\right),
            \nonumber \\  
  \Delta^{\lambda,\mu}_{w,\nu}(P_{\pm}) & = & \frac{\lambda}{\nu}\, P_{\pm}
            \otimes e^{\mu w P_3} + \frac{\mu}{\nu}\, e^{-\lambda w P_3} \otimes
            P_{\pm}, \nonumber \\
  \epsilon_{w,\nu}(X) & = & 0, \qquad X \in \{J_3, J_{\pm}, P_3, P_{\pm}\},
            \nonumber \\ 
  S^{\mu}_{w,\nu}(J_3) & = & - J_3, \qquad S^{\mu}_{w,\nu}(P_3) = -
            \frac{\mu}{\nu} P_3, \nonumber \\
  S^{\mu}_{w,\nu}(J_{\pm}) & = & - \left(J_{\pm} \pm 2\mu w P_{\pm}\right),
            \qquad S^{\mu}_{w,\nu}(P_{\pm}) = - \frac{\mu}{\nu} P_{\pm}, \nonumber \\
  {\cal R}^{\lambda,\mu}_w & = & \exp\left\{2w \left(\lambda P_3 \otimes J_3 +
            \mu J_3 \otimes P_3\right)\right\} \exp\left\{B^{\lambda,\mu}
            \arcsinh\left(2w A^{\lambda,\mu}\right) \Big/ \left(w
            A^{\lambda,\mu}\right)\right\} \nonumber \\
  & & \mbox{} \times \left(1 + 4 w^2 \left(A^{\lambda,\mu}\right)^2\right)^{-1/2}, 
\end{eqnarray}
where
\begin{eqnarray}
  A^{\lambda,\mu} & \equiv & w Q^{\lambda}_+ \otimes Q^{\mu}_-, \nonumber \\
  B^{\lambda,\mu} & \equiv & w \left( L^{\lambda}_+ \otimes Q^{\mu}_- +
             Q^{\lambda}_+ \otimes L^{\mu}_-\right) - w^2 \left(2 Q^{\lambda}_+
             \otimes Q^{\mu}_- + Q^{\lambda}_+ \otimes J_3 Q^{\mu}_- - J_3
             Q^{\lambda}_+ \otimes Q^{\mu}_-\right), \nonumber \\
  L^{\lambda}_{\pm} & \equiv & e^{\pm\lambda w P_3} J_{\pm}, \qquad
  Q^{\lambda}_{\pm} \equiv \lambda e^{\pm\lambda w P_3} P_{\pm}. 
\end{eqnarray}
\par
%
%
The quantum Euclidean algebra $U_w(e(3))$ admits the $4 \times 4$ matrix
representation
\begin{eqnarray}
  D(J_3) & = & -i e_{12} + i e_{21}, \qquad D(J_{\pm}) = \mp e_{13} - i e_{23} \pm
            e_{31} + i e_{32}, \nonumber \\
  D(P_3) & = & e_{34}, \qquad D(P_{\pm}) = e_{14} \pm i e_{24},
\end{eqnarray}
where $e_{ij}$ denotes the matrix with entry~1 in row~$i$ and column~$j$, and
zeros everywhere else. In such a representation, the coloured universal $\cal
R$-matrix is represented by the $16 \times 16$ matrix
\begin{eqnarray}
  R^{\lambda,\mu}_w & \equiv & (D \otimes D)\left({\cal R}^{\lambda,\mu}_w\right)
          \nonumber \\
  & = &\left(\begin{array}{cccc}
          1_4 & -2i\mu w D(P_3) & -2\mu w D(P_-) & 2\lambda w D(J_-) \\[0.1cm]
          2i\mu w D(P_3) & 1_4 & -2i\mu w D(P_-) & 2i\lambda w D(J_-) \\[0.1cm]
          2\mu w D(P_-) & 2i\mu w D(P_-) & 1_4 & 2\lambda w D(J_3) \\[0.1cm]
          0_4 & 0_4 & 0_4 & 1_4
  \end{array}\right),   
\end{eqnarray}
which is a new solution of the coloured YBE.\par
%
%
{\it Remark.\/} Similar results can be obtained for the two-dimensional quantum
Euclidean algebra $U_w(e(2))$~\cite{celeghini90}, but in such a case no coloured
universal $\cal R$-matrix is known.\par
%
%
\subsection{Null-plane $D$-dimensional quantum Poincar\'e algebras
$U_z(iso(D-1,1))$}   \label{subsec:poincare}
As a last example, we consider the  null-plane deformations $U_z(iso(D-1,1))$ of
the Poincar\'e algebras in $D=2$~\cite{ballesteros95a},
$D=3$~\cite{ballesteros95b}, and $D=4$ dimensions~\cite{ballesteros95c}. Since
the results look quite similar for different $D$~values, we only list here those
for~$D=4$.\par
%
%
The quantum algebra $U_z(iso(3,1))$ is generated by ten operators $K_3$, $J_3$,
$P_+$, $P_-$, $P_1$, $P_2$, $E_1$, $E_2$, $F_1$,~$F_2$, which, in the $z\to 0$
limit, go over into the following combinations of Poincar\'e algebra generators in
the usual physical basis, $J_i$ (rotations), $K_i$ (boosts), and $P_{\mu}$
(translations), where $i=1$, 2,~3, and $\mu=0$, 1, 2,~3: $P_{\pm} = (P_0 \pm
P_3)/2$, $E_1 = (K_1 + J_2)/2$, $E_2 = (K_2 - J_1)/2$, $F_1 = (K_1 - J_2)/2$,
$F_2 = (K_2 + J_1)/2$. Their nonvanishing commutators are given by
\begin{eqnarray}
  \left[K_3, P_+\right] & = & \frac{e^{2z P_+} -1}{2z}, \qquad \left[K_3, P_-\right] =
            - P_- - z P_1^2 - z P_2^2, \nonumber \\
  \left[K_3, E_i\right] & = & E_i e^{2z P_+}, \qquad \left[K_3, F_i\right] =
            - F_i - 2z K_3 P_i, \nonumber \\
  \left[J_3, P_i\right] & = & - \epsilon_{ij3} P_j, \qquad  \left[J_3, E_i\right] = 
            - \epsilon_{ij3} E_j, \qquad  \left[J_3, F_i\right] = - \epsilon_{ij3}
            F_j, \nonumber \\
  \left[E_i, P_j\right] & = & \delta_{ij} \frac{e^{2z P_+} -1}{2z}, \qquad \left[F_i,
            P_j\right] = \delta_{ij} \left(P_- + z P_1^2 + z P_2^2\right), \nonumber \\
  \left[E_i, F_j\right] & = & \delta_{ij} K_3 + \epsilon_{ij3} J_3 e^{2z P_+}, \qquad
            \left[P_+, F_i\right] = - P_i, \nonumber \\
  \left[F_1, F_2\right] & = & 2z \left(P_1 F_2 - P_2 F_1\right), \qquad
            \left[P_-, E_i\right] = - P_i,   \label{eq:poincare-def}
\end{eqnarray}
where $i$,~$j$ run over 1,~2. Here we assume $k = \R$, and $z \in {\cal Q} = \R
\! \setminus \! \{0\}$.\par
%
%
The algebra defining relations~(\ref{eq:poincare-def}) are left invariant under the
transformations
\begin{eqnarray}
  \sigma^{\nu}(K_3) & = & K_3, \qquad \sigma^{\nu}(J_3) = J_3, \qquad
            \sigma^{\nu}(P_+) = \nu_1 \nu_2 P_+, \qquad \sigma^{\nu}(P_-) = \nu_1 
            \nu_2^{-1} P_-, \nonumber \\
  \sigma^{\nu}(P_i) & = & \nu_1 P_i, \qquad \sigma^{\nu}(E_i) = \nu_2 E_i, \qquad
            \sigma^{\nu}(F_i) = \nu_2^{-1} F_i, 
\end{eqnarray}
where $\nu \equiv (\nu_1, \nu_2) \in {\cal C} = (\R \! \setminus \! \{0\}) \times (\R
\! \setminus \! \{0\})$, provided $z$ is changed into $z^{\nu} = \nu_1 \nu_2 z$. The
corresponding colour group is ${\cal G} = Gl(1,\R) \otimes Gl(1,\R)$.\par
%
%
The coloured maps and universal $\cal R$-matrix are found to be given by
\begin{eqnarray}
  \Delta^{\lambda,\mu}_{z,\nu}(J_3) & = & J_3 \otimes 1 + 1 \otimes J_3, \qquad
            \Delta^{\lambda,\mu}_{z,\nu}(P_+) =
            \frac{\lambda_1\lambda_2}{\nu_1\nu_2}\,
            P_+ \otimes 1 + \frac{\mu_1\mu_2}{\nu_1\nu_2}\, 1 \otimes P_+,
            \nonumber \\
  \Delta^{\lambda,\mu}_{z,\nu}(P_-) & = & \frac{\lambda_1\nu_2}{\lambda_2\nu_1}
            \, P_- \otimes e^{2\mu_1\mu_2 z P_+} + \frac{\mu_1\nu_2}{\mu_2\nu_1}\,
           1 \otimes P_-, \nonumber \\
  \Delta^{\lambda,\mu}_{z,\nu}(P_i) & = & \frac{\lambda_1}{\nu_1}\, P_i \otimes
            e^{2\mu_1\mu_2 z P_+} + \frac{\mu_1}{\nu_1}\, 1 \otimes P_i, \qquad
            \Delta^{\lambda,\mu}_{z,\nu}(E_i) =\frac{\lambda_2}{\nu_2}\, E_i
            \otimes 1 + \frac{\mu_2}{\nu_2}\, 1 \otimes E_i, \nonumber \\
  \Delta^{\lambda,\mu}_{z,\nu}(F_1) & = & \frac{\nu_2}{\lambda_2}\, F_1 \otimes
            e^{2\mu_1\mu_2 z P_+} + \frac{\nu_2}{\mu_2}\, 1 \otimes F_1 -
            \frac{2\lambda_1\mu_2\nu_2}{\lambda_2} z P_- \otimes E_1\, 
            e^{2\mu_1\mu_2 z P_+} \nonumber \\
  & & \mbox{} - 2\lambda_1\nu_2 z P_2 \otimes J_3\, e^{2\mu_1\mu_2 z P_+},
            \nonumber \\   
  \Delta^{\lambda,\mu}_{z,\nu}(F_2) & = & \frac{\nu_2}{\lambda_2}\, F_2 \otimes
            e^{2\mu_1\mu_2 z P_+} + \frac{\nu_2}{\mu_2}\, 1 \otimes F_2 -
            \frac{2\lambda_1\mu_2\nu_2}{\lambda_2} z P_- \otimes E_2\, 
            e^{2\mu_1\mu_2 z P_+} \nonumber \\
  & & \mbox{} + 2\lambda_1\nu_2 z P_1 \otimes J_3\, e^{2\mu_1\mu_2 z P_+},
            \nonumber \\
  \Delta^{\lambda,\mu}_{z,\nu}(K_3) & = & K_3 \otimes e^{2\mu_1\mu_2 z P_+} +
            1 \otimes K_3 - 2\lambda_1\mu_2 z P_1 \otimes E_1\, 
            e^{2\mu_1\mu_2 z P_+} \nonumber \\
  & & \mbox{} - 2\lambda_1\mu_2 z P_2 \otimes E_2\, e^{2\mu_1\mu_2 z P_+},
            \nonumber \\
  \epsilon_{z,\nu}(X) & = & 0, \qquad X \in \{K_3, J_3, P_{\pm}, P_i, E_i, F_i\},
            \nonumber \\
  S^{\mu}_{z,\nu}(J_3) & = & - J_3, \qquad S^{\mu}_{z,\nu}(P_+) = -
            \frac{\mu_1\mu_2}{\nu_1\nu_2} P_+, \qquad S^{\mu}_{z,\nu}(P_-) = -
            \frac{\mu_1\nu_2}{\mu_2\nu_1} P_- \, e^{-2\mu_1\mu_2 z P_+},
            \nonumber \\
  S^{\mu}_{z,\nu}(P_i) & = & - \frac{\mu_1}{\nu_1} P_i\, e^{-2\mu_1\mu_2 z P_+},
            \qquad \qquad S^{\mu}_{z,\nu}(E_i) = - \frac{\mu_2}{\nu_2} E_i,
            \nonumber \\
  S^{\mu}_{z,\nu}(F_1) & = & - \frac{\nu_2}{\mu_2} \left(F_1 + 2\mu_1\mu_2 z P_-
            E_1 + 2\mu_1\mu_2 z P_2 J_3\right) e^{-2\mu_1\mu_2 z P_+},
            \nonumber \\
  S^{\mu}_{z,\nu}(F_2) & = & - \frac{\nu_2}{\mu_2} \left(F_2 + 2\mu_1\mu_2 z P_-
            E_2 - 2\mu_1\mu_2 z P_1 J_3\right) e^{-2\mu_1\mu_2 z P_+},
            \nonumber \\
  S^{\mu}_{z,\nu}(K_3) & = & - \left(K_3 + 2\mu_1\mu_2 z P_1 E_1 + 2\mu_1\mu_2
            z P_2 E_2\right) e^{-2\mu_1\mu_2 z P_+}, \nonumber \\
  {\cal R}^{\lambda,\mu}_z & = & \exp\left\{2\lambda_2\mu_1 z E_2 \otimes
            P_2\right\} \exp\left\{2\lambda_2\mu_1 z E_1 \otimes P_1\right\}
            \exp\left\{-2\lambda_1\lambda_2 z P_+ \otimes K_3 \right\}
            \nonumber \\
  & & \mbox{} \times \exp\left\{2\mu_1\mu_2 z K_3 \otimes P_+\right\}
            \exp\left\{-2\lambda_1\mu_2 z P_1 \otimes E_1 \right\} \nonumber \\
  & & \mbox{} \times \exp\left\{-2\lambda_1\mu_2 z P_2 \otimes E_2\right\}.
\end{eqnarray}
\par
%
%
The quantum Poincar\'e algebra $U_z(iso(3,1))$ admits the $5 \times 5$ matrix
representation
\begin{eqnarray}
  D(K_3) & = & e_{14} + e_{41}, \qquad D(J_3) = e_{23} - e_{32}, \qquad D(P_+) =
           \case{1}{2} \left(e_{10} + e_{40}\right), \nonumber \\
  D(P_-) & = & e_{10} - e_{40}, \qquad D(P_1) = e_{20}, \qquad D(P_2) = e_{30},
           \nonumber \\
  D(E_1) & = & \case{1}{2} \left(e_{12} + e_{21} - e_{24} + e_{42}\right), \qquad
           D(E_2) = \case{1}{2} \left(e_{13} + e_{31} - e_{34} + e_{43}\right), 
           \nonumber \\
  D(F_1) & = & e_{12} + e_{21} + e_{24} - e_{42}, \qquad D(F_2) = e_{13} + e_{31} +
           e_{34} - e_{43}, 
\end{eqnarray}
where rows and columns are labelled by 0, 1, 2, 3,~4, and $e_{ij}$ has the same
meaning as in Subsec.~\ref{subsec:euclide}. In such a representation, the coloured
universal $\cal R$-matrix gives rise to the following new $25 \times 25$ matrix
solution of the coloured YBE:
\begin{eqnarray}
  R^{\lambda,\mu}_z & \equiv & (D \otimes D)\left({\cal R}^{\lambda,\mu}_z\right)
          \nonumber \\
  & = &\left(\begin{array}{ccccc}
          \ss 1_5 & \ss 0_5 & \ss 0_5 & \ss 0_5 & \ss 0_5 \\
          \ss -\lambda_1\lambda_2 z D(K_3) & \ss 1_5 & \ss \lambda_2\mu_1 z
                D(P_1) & \ss \lambda_2\mu_1 z D(P_2) & \ss 2\mu_1\mu_2 z D(P_+) \\
          \ss -2\lambda_1\mu_2 z D(E_1) & \ss \lambda_2\mu_1 z D(P_1) & \ss 1_5
                & \ss 0_5 & \ss -\lambda_2\mu_1 z D(P_1) \\
          \ss -2\lambda_1\mu_2 z D(E_2) & \ss \lambda_2\mu_1 z D(P_2) & \ss 0_5
                & \ss 1_5 & \ss -\lambda_2\mu_1 z D(P_2) \\
          \ss -\lambda_1\lambda_2 z D(K_3) & \ss 2\mu_1\mu_2 z D(P_+) & \ss
                \lambda_2\mu_1 z D(P_1) & \ss \lambda_2\mu_1 z D(P_2) & \ss 1_5 
  \end{array}\right).
\end{eqnarray}
\par
%
%
{\it Remark.\/} The $\kappa$-deformations $U_{\kappa}(iso(D-1,1))$ of the
$D$-dimensional Poincar\'e algebras~\cite{lukierski} can be transformed into
coloured Hopf algebras along the same lines as $U_z(iso(D-1,1))$, but in such a
case no coloured universal $\cal R$-matrix is known.\par
%
%
\section{CONCLUSION}   \label{sec:conclusion}
In the present paper, we did introduce some new algebraic structures, termed
coloured Hopf algebras, by combining the coalgebra structures and antipodes of a
standard Hopf algebra set with the transformations of an algebra isomorphism
group, called colour group. We did show that various classes of Hopf algebras,
such as almost cocommutative, coboundary, quasitriangular, and triangular ones,
can be extended into corresponding coloured structures, and that coloured
quasitriangular Hopf algebras, in particular, are characterized by the existence of a
coloured universal $\cal R$-matrix, satisfying the coloured YBE.\par
%
%
Finally, we did apply the new concepts to QUEA's of both semisimple and
nonsemisimple Lie algebras, and did prove by means of examples that the colour
group may be chosen as a finite or infinite, abelian or nonabelian group. Through
such constructions, we did demonstrate that the coloured Hopf algebras defined
here significantly generalize those previously introduced by Ohtsuki~\cite{ohtsuki},
because the latter are restricted to abelian colour groups, in which case they
reduce to substructures of the former.\par
%
%
It is worth noting that some of the matrix representations of coloured universal
$\cal R$-matrices constructed in the present paper, as well as those that would be
obtained in higher-dimensional representations, provide new solutions of the
coloured YBE, which might be of interest in the context of integrable models.\par
%
%
It is also important to stress that the applicability of the coloured Hopf algebra
new concept is not confined to QUEA's of Lie algebras. As we plan to show
elsewhere, QUEA's of Lie superalgebras may also provide a suitable starting point
for constructing coloured Hopf algebras.\par
%
%
Other types of Hopf algebras might be used as well, such as those arising in the
FRT~formalism. The resulting coloured algebraic structures would significantly
differ from those previously constructed by Kundu and Basu-Mallick~\cite{basu94,
kundu94, basu95}, since the latter have the same coalgebra structure as the original
Hopf algebras, whereas for the former it is the algebra structure that would be left
unchanged. Further investigation of possible relationships between both types of
coloured algebraic structures would be highly desirable.\par
%
%
In the examples considered in the present paper, no effort has been made to
determine the maximal colour group --- hence the maximal coloured Hopf structure
--- compatible with a given Hopf algebra set. Similarly, the restrictions on the
colour parameters imposed by considering a given real form of a complex Hopf
algebra have not been systematically investigated. Solving such problems might be
interesting topics for future study.\par
%
%
Generalizing to coloured algebraic structures the duality relationship between 
pairs of Hopf algebras $U_q(g)$ and~$G_q$, as highlighted in the universal $\cal
T$-matrix formalism~\cite{reshetikhin}, might also be a promising
direction for future investigation.\par
%
%
\section*{ACKNOWLEDGMENT}
The author would like to thank C.~Daskaloyannis and L.~Hlavat\'y for some valuable
comments.\par
\newpage
%
%


\begin{thebibliography}{99}
%
\bibitem{yang} C. N. Yang, Phys. Rev. Lett. {\bf 19}, 1312 (1967); R. J. Baxter, {\it
Exactly Solved Models in Statistical Mechanics\/} (Academic, New York, 1982).
%
\bibitem{kulish} P. Kulish and E. K. Sklyanin, in {\it Integrable Quantum Field
Theories, Tvarminne, 1981\/}, edited by J. Hietarinta and C. Montonen, Lecture
Notes in Physics, Vol.~151 (Springer, Berlin, 1982) p.~61; L. Faddeev, in {\it
Integrable Systems\/}, edited by M.-L. Ge and X. C. Song (World Scientific,
Singapore, 1989) p.~23; L. Alvarez-Gaum\'e, C. G\'omez, and G. Sierra, Phys. Lett. B
{\bf 220}, 142 (1990); H. J. de Vega, Int. J. Mod. Phys. A {\bf 4}, 2371 (1989); E.
Witten, Nucl. Phys. B {\bf 330}, 285 (1990).
%
\bibitem{jones} V. F. R. Jones, Bull. Amer. Math. Soc. {\bf 12}, 103 (1985); Y.
Akutsu and M. Wadati, J. Phys. Soc. Jpn. {\bf 56}, 839 (1987); V. G. Turaev, Invent.
Math. {\bf 92}, 527 (1988).
%
\bibitem{sklyanin} E. K. Sklyanin, Funct. Anal. Appl. {\bf 16}, 262 (1982); P. P. Kulish 
and N. Y. Reshetikhin, J. Sov. Math. {\bf 23}, 2435 (1983).  
%
\bibitem{faddeev} L. Faddeev, N. Reshetikhin, and L. Takhtajan, in {\it Algebraic
Analysis\/}, Vol.~1, edited by M. Kashiwara and T. Kawai (Academic, New York,
1988) p.~129; in {\it Braid Group, Knot Theory and Statistical Mechanics\/}, edited
by C. N. Yang and M. L. Ge (World Scientific, Singapore, 1989) p.~97.
%
\bibitem{drinfeld} V. G. Drinfeld, in {\it Proc. Int. Congress of
Mathematicians (Berkeley, CA, 1986)\/}, edited by A. M. Gleason (AMS, Providence,
RI, 1987), p. 798; M. Jimbo, Lett. Math. Phys. {\bf 10}, 63 (1985); {\bf 11}, 247
(1986). 
%
\bibitem{celeghini90} E. Celeghini, R. Giachetti, E. Sorace, and M. Tarlini, J. Math.
Phys. {\bf 31}, 2548 (1990); {\bf 32}, 1155 (1991).
%
\bibitem{celeghini91} E. Celeghini, R. Giachetti, E. Sorace, and M. Tarlini, J. Math.
Phys. {\bf 32}, 1159 (1991).
%
\bibitem{schirrmacher} A. Schirrmacher, J. Wess, and B. Zumino, Z. Phys. C {\bf 49},
317 (1991); A. Schirrmacher, {\it ibid.\/} {\bf 50}, 321 (1991).
%
\bibitem{dobrev} V. K. Dobrev, J. Math. Phys. {\bf 33}, 3419 (1992); J. Geom. Phys.
{\bf 11}, 367 (1993).
%
\bibitem{burroughs} N. Burroughs, Commun. Math. Phys. {\bf 133}, 91 (1990).
%
\bibitem{bazhanov} V. V. Bazhanov and Yu. G. Stroganov, Theor. Math. Phys. {\bf 62},
253 (1985); R. J. Baxter, J. H. Perk, and H. Au-Yang, Phys. Lett. A {\bf 128}, 138
(1988); B. S. Shastry, J. Stat. Phys. {\bf 50}, 57 (1988).
%
\bibitem{hlavaty} L. Hlavat\'y, J. Phys. A {\bf 20}, 1661 (1987); J. Hietarinta and C.
Viallet, ``On the parametrization of solutions of the Yang-Baxter equations,''
Universit\'e Pierre et Marie Curie preprint PAR-LPTHE 94-25, q-alg/9504028
(1995).
%
\bibitem{footnote1} The coloured YBE considered here should not be confused with
the colour YBE~\cite{mcanally} that arises when extending the graded
YBE~\cite{chaichian} to more general gradings than that determined by
$\Z_2$~\cite{rittenberg}. Similarly, the coloured Hopf algebras introduced in
Sec.~\ref{sec:colHopf}, as well as those considered in Refs.~\cite{bonatsos,
ohtsuki}, are distinct from the Hopf colour algebras~\cite{mcanally}, generalizing
Hopf superalgebras~\cite{chaichian} to such more general gradings.
%
\bibitem{akutsu} Y. Akutsu and T. Deguchi, Phys. Rev. Lett. {\bf 67}, 777 (1991);
Mod. Phys. Lett. A {\bf 7}, 767 (1992).
%
\bibitem{ge} M.-L. Ge, X.-F. Liu, and C.-P. Sun, Phys. Lett. A {\bf 155}, 137 (1991);
M.-L. Ge, C.-P. Sun, and K. Xue, Int. J. Mod. Phys. A {\bf 7}, 6609 (1992).
%
\bibitem{burdik92a} \v C. Burd\'\i k and P. Hellinger, J. Phys. A {\bf 25}, L1023
(1992).
%
\bibitem{gomez} C. G\'omez and G. Sierra, J. Math. Phys. {\bf 34},  2119 (1993).
%
\bibitem{basu94} B. Basu-Mallick, Mod. Phys. Lett. A {\bf 9}, 2733 (1994).
%
\bibitem{kundu94} A. Kundu and B. Basu-Mallick, J. Phys. A {\bf 27}, 3091 (1994). 
%
\bibitem{kundu92} A. Kundu and B. Basu-Mallick, J. Phys. A {\bf 25}, 6307 (1992).
%
\bibitem{bonatsos} D. Bonatsos, C. Daskaloyannis, P. Kolokotronis, A. Ludu, and C.
Quesne, J. Math. Phys. {\bf 38}, 369 (1997); D. Bonatsos, P. Kolokotronis, C.
Daskaloyannis, A. Ludu, and C. Quesne, Czech. J. Phys. {\bf 46}, 1189 (1996).
%
\bibitem{basu95} B. Basu-Mallick, Int. J. Mod. Phys. A {\bf 10}, 2851 (1995).
%
\bibitem{ohtsuki} T. Ohtsuki, J. Knot Theor. Its Rami. {\bf 2}, 211 (1993).
%
\bibitem{majid} S. Majid, Int. J. Mod. Phys. A {\bf 5}, 1 (1990); V. Chari and A.
Pressley, {\it A Guide to Quantum Groups} (Cambridge U.P., Cambridge, 1994).
%
\bibitem{kirillov} A. N. Kirillov and N. Yu. Reshetikhin, in {\it Infinite-dimensional
Lie Algebras and Groups\/}, edited by V. G. Kac (World Scientific, Singapore, 1989)
p.~285.
%
\bibitem{burdik92b} \v C. Burd\'\i k and P. Hellinger, J. Phys. A {\bf 25}, L629
(1992).
%
\bibitem{ohn} C. Ohn, Lett. Math. Phys. {\bf 25}, 85 (1992).
%
\bibitem{agha} A. Aghamohammadi, M. Khorrami, and A. Shariati, J. Phys. A {\bf 28},
L225 (1995).
%
\bibitem{shariati} A. Shariati, A. Aghamohammadi, and M. Khorrami, Mod. Phys. Lett.
A {\bf 11}, 187 (1996).
%
\bibitem{ballesteros96a} A. Ballesteros and F. J. Herranz, J. Phys. A {\bf 29}, L311
(1996).
%
\bibitem{footnote2} Note that Eq.~(24) in Ref.~\cite{ballesteros96a} is not correct
and is replaced by Eq.~(\ref{eq:ns-sl(2)-rep}) of the present paper.
%
\bibitem{abdesselam} B. Abdesselam, A. Chakrabarti, and R. Chakrabarti,
``Irreducible representations of Jordanian quantum algebra $U_h(sl(2))$ via a
nonlinear map,'' Ecole Polytechnique preprint CPTH-S455.0696, q-alg/9606014
(1996).
%
\bibitem{ballesteros96b} A. Ballesteros and F. J. Herranz, J. Phys. A {\bf 29}, 4307
(1996).
%
\bibitem{lyakhovsky} V. Lyakhovsky and A. Mudrov, J. Phys. A {\bf 25}, L1139 (1992).
%
\bibitem{ballesteros95a} A. Ballesteros, E. Celeghini, F. J. Herranz, M. A. del Olmo,
and M. Santander, J. Phys. A {\bf 28}, 3129 (1995); A. Ballesteros, F. J. Herranz, M. A.
del Olmo, C. M. Pere\~na, and M. Santander, {\it ibid.\/} A {\bf 28}, 7113 (1995); M.
Khorrami, A. Shariati, M. R. Abolhassani, and A. Aghamohammadi, Mod. Phys. Lett. A
{\bf 10}, 873 (1995). 
%
\bibitem{ballesteros95b} A. Ballesteros, F. J. Herranz, M. A. del Olmo, and M.
Santander, J. Phys. A {\bf 28}, 941 (1995); A. Ballesteros and F. J. Herranz,
``$(2+1)$ null-plane quantum Poincar\'e group from a factorized universal 
$R$-matrix,'' Universidad de Burgos preprint q-alg/9605031 (1996); A. Shariati, A.
Aghamohammadi, and M. Khorrami,  Mod. Phys. Lett. A {\bf 11}, 187 (1996). 
%
\bibitem{ballesteros95c} A. Ballesteros, F. J. Herranz, M. A. del Olmo, and M.
Santander, Phys. Lett. B {\bf 351}, 137 (1995); A. Ballesteros, F. J. Herranz, and C. M.
Pere\~na, {\em ibid.} B {\bf 391}, 71 (1997).
%
\bibitem{lukierski} J. Lukierski, A. Nowicki, and H. Ruegg, Phys. Lett. B {\bf 293},
344 (1992); J. Lukierski and H. Ruegg, {\it ibid.\/} B {\bf 329}, 189 (1994).
%
\bibitem{reshetikhin} N. Reshetikhin and M. Semenov-Tian-Shansky, J. Geom. Phys.
{\bf 5}, 533 (1988); C. Fronsdal and A. Galindo, Lett. Math. Phys. {\bf 27}, 59
(1993); F. Bonechi, E. Celeghini, R. Giachetti, C. M. Pere\~na, E. Sorace, and M.
Tarlini, J. Phys. A {\bf 27}, 1307 (1994).
%
\bibitem{mcanally} D. S. McAnally, in {\it Proc. Yamada Conf. XL, XX Int. Coll. on
Group Theoretical Methods in Physics, Toyonaka, Japan, July 4--9, 1994\/}, edited
by A. Arima, T. Eguchi, and N. Nakanishi (World Scientific, Singapore, 1995) p.~339.
%
\bibitem{chaichian} M. Chaichian and P. Kulish, Phys. Lett. B {\bf 234}, 72 (1990).
%
\bibitem{rittenberg} V. Rittenberg and D. Wyler, Nucl. Phys. B {\bf 139}, 189
(1978); J. Math. Phys. {\bf 19}, 2193 (1978); J. Lukierski and V. Rittenberg, Phys.
Rev. D {\bf 18}, 385 (1978); M. Scheunert, J. Math. Phys. {\bf 20}, 712 (1979).

\end{thebibliography}
\end{document}